\documentclass[a4paper,reqno,9pt,oneside]{amsart}
\usepackage{amsmath,geometry}
\usepackage{color}

\numberwithin{equation}{section}

\newtheorem{thm}{Theorem}[section]

\newtheorem{pro}[thm]{Proposition}
\newtheorem{lm}[thm]{Lemma}
\newtheorem{oss}[thm]{Remark}

\begin{document}

\date{}

\title[Singular limits in higher order Lioville-type equations]{Singular limits in higher order Lioville-type equations}

\author{Fabrizio Morlando}
\address{Fabrizio Morlando, Dipartimento di Matematica e Fisica, Universit\`a degli Studi Roma Tre,\,Largo S. Leonardo Murialdo 1, 00146 Roma, Italy}
\email{morlando@mat.uniroma3.it}
\begin{abstract}
In this paper we consider the following higher order boundary value problem
\medskip
$$
\left\{ \begin{array}{ll}
          (-\Delta)^{m} u=\rho^{2m} V(x) e^{u}& \mbox{in}\ \Omega\\
          B_{j}u=0, |j|\leq m-1& \mbox{on}\ \partial\Omega,\
        \end{array} \right.
\vspace{0,1cm}$$
where $\Omega$ is a smooth bounded domain in $\mathbb{R}^{2m}$ with $m\in\mathbb{N}$, $V(x)\neq0$ is a smooth function positive somewhere in $\Omega$ and $\rho$ is a positive small parameter. Here, the operator $B_{j}$ stands for either Navier or Dirichlet boundary conditions. We find sufficient conditions under which, as $\rho$ approaches $0$, there exists an explicit class of solutions which admit a concentration behavior with a prescribed bubble profile around some given $k$-points in $\Omega$, for any given integer $k$. These are the so-called singular limits. The candidate $k$-points of concentration must be critical points of a suitable finite dimensional functional explicitly defined in terms of the potential $V$ and the higher order Green's function with respect to the imposed boundary conditions.
\end{abstract}

\maketitle

\section{Introduction}

Let $\Omega \subseteq\mathbb{R}^{2m}$ be a smooth bounded domain with $m\in\mathbb{N}$. This paper deals with the existence, the qualitative properties and the asymptotic behavior of nontrivial solutions $u: \Omega\rightarrow\mathbb{R}$ to the following higher order boundary value problem
\begin{equation}\label{E1}
\left\{ \begin{array}{ll}
          (-\Delta)^{m} u=\rho^{2m} V(x) e^{u}& \mbox{in}\ \Omega\\
          B_{j}u=0, |j|\leq m-1& \mbox{on}\ \partial\Omega,\
        \end{array} \right.
\vspace{0,2cm}\end{equation}
where we prescribe the boundary conditions to be either Navier $B_{j}u:=(-\Delta)^{j}u$ or Dirichlet $B_{j}u:=\partial^{j}_{\nu}u$ with $\nu=\nu(x)$ the unit outer normal, $V(x)\neq0$ a given smooth potential and $\rho\in\mathbb{R}^{+}$ a small parameter which tends to zero from above. The leading part $(-\Delta)^{m}$  is, usually, called \emph{polyharmonic operator} on $\mathbb{R}^{2m}$, it is simply obtained as the composition $m$-times of $-\Delta$. In the sequel, we denote by $\mathcal{H}_{m}(\Omega)$ the linear space of polyharmonic functions of \emph{order} $m$ in $\Omega$. We recall that \eqref{E1} corresponds to a standard case of \emph{uniform singular convergence}, in the sense that the associated nonlinear coefficient $\rho^{2m}k(x)$ goes to zero uniformly in $\overline{\Omega}$ as $\rho\downarrow0$. In the last decades a lot of work has been done in the study of the asymptotic behavior of solutions in the limit, the so-called \emph{singular limits}.
\\Higher order elliptic equations involving exponential nonlinearities appear naturally in conformal geometry. On a compact surface $(M,g)$ with a Riemannian metric $g$, a natural curvature invariant associated with the Laplace-Beltrami operator is the Gaussian curvature $K=K_{g}$. Under the conformal change of metric $\widetilde{g}=e^{2u}g$, we have
$$-\Delta_{g}+K=K_{u}e^{2u}\ \mbox{on}\ M$$ where $K_{u}$ denotes the Gaussian curvature of $(M,\widetilde{g})$. For compact manifolds of general dimension $N$, when $N$ is even, the existence of a $N$-th order operator $\mathbf{P}_{g}^{N}$ conformally covariant of bi-degree $(0,N)$ was verified in \cite{gjms}. This operator is usually referred to
as the \emph{GJMS} operator and its leading term is the polyharmonic operator. However, this operator is only explicitly known on the standard Euclidean space $\mathbb{R}^{N}$ where it reduces to its leading term and on the standard sphere $\mathbb{S}^{N}$ through the Beckner-Branson formula \cite{beck,bran3}.
Attached to this operator, there is a natural concept of curvature, namely the \emph{$Q_{g}^{N}$-curvature}, firstly introduced by Branson in \cite{bran}.
If $N=2m$ is even, we have after a conformal change of metric $\tilde{g}=e^{2u}g$ that
\begin{equation}
\mathbf{P}_{\tilde{g}}^{N}=e^{-Nu}\mathbf{P}_{g}^{N},
\end{equation}
and
\begin{equation}\label{eq}
\mathbf{P}_{g}^{N}u+\mathbf{Q}_{g}^{N}=\mathbf{Q}_{\tilde{g}}^{N}e^{Nu}.
\end{equation}
We notice that \eqref{eq} is a generalized version of Gauss's identity in dimension $2$.\\
\\One can think the \emph{GJMS}-operator as a higher order analog of the Laplace-Beltrami operator defined on a compact $2$-manifold, and the $Q$-curvature can be thought of as a higher order analog of the Gaussian curvature $K_{g}$, in fact, in dimension $2$ we simply have $\mathbf{P}_{g}^{2}=-\Delta_{g}$ and $\mathbf{Q}_{g}^{2}=K_{g}$. We refer to the research monograph \cite{014} for related topics and recent developments.
Therefore, in the spirit of the classical Poincar\'{e} uniformization theorem for closed surfaces, a similar question arises:
\emph{for every closed $N$-dimensional manifold $(M,g)$ with $N$ even, does there exist a metric $\tilde{g}$, conformally equivalent to $g$, for which the corresponding Q-curvature $\mathbf{Q}_{\tilde{g}}^{N}$ is constant?}
\\In view of the transformation law $\tilde{g}=e^{2u}g$, the problem reduces to solve the following $N$-th order nonlinear elliptic equation:
\begin{equation}\label{E2}
\mathbf{P}_{g}^{N}u+\mathbf{Q}_{g}^{N}=\mathbf{\overline{Q}}e^{Nu}\ \mbox{in}\ M,
\end{equation}
where $\mathbf{\overline{Q}}$ is a real constant. The above equation has been extensively studied during the last decades. In (\cite{ndi}, Theorem 1.1), generalizing a result in a seminal paper by Djadli and Malchiodi in the fourth dimensional case \cite{0120}, the author solves completely this problem under generic assumptions: $$\ker\mathbf{P}_{g}^{N}=\{\mbox{constant}\},\int_{M}\mathbf{Q}_{g}^{N}dv_{g}\neq k(N-1)!\omega_{N}\ \mbox{for}\ k\in\mathbb{N}^{\ast}.$$
Here we are interested in the special case where the manifold is the Euclidean space $\mathbb{R}^{2m}$ endowed with the Euclidean metric $g_{\mathbb{R}^{2m}}$. If $N=2m$, $\mathbf{P}_{g_{\mathbb{R}^{2m}}}^{2m}$ are pointwise operators and by stereographic projection on $\mathbb{R}^{2m}$, they simply become $(-\Delta)^{m}$, whereas, if $N$ is odd, $\mathbf{P}_{g_{\mathbb{R}^{N}}}^{N}$ are nonlocal pseudo differential operators and by stereographic projection on $\mathbb{R}^{N}$, they simply become $(-\Delta)^{\frac{N}{2}}$. In both cases, by the flatness of $\mathbb{R}^{N}$, we have $\mathbf{Q}_{g_{\mathbb{R}^{N}}}^{N}\equiv0$. In the case $N$ even, \eqref{E2} reduces to
\begin{equation}\label{E3}
(-\Delta)^{m}u=\mathbf{\overline{Q}}e^{2mu}\ \mbox{in}\ \mathbb{R}^{2m},
\end{equation}
which gives rise to a conformal metric $\tilde{g}=e^{2u}g_{\mathbb{R}^{2m}}$ whose $Q$-curvature is given by the constant $\mathbf{\overline{Q}}$. Clearly we can set \eqref{E3} in a regular bounded domain $\Omega\subseteq\mathbb{R}^{2m}$ but, since we are no more dealing with a manifold without a boundary, we need to impose some boundary conditions on $\partial\Omega$ in order to have a well-posed problem.
The main analytic difficult in working with nonlinear higher order elliptic equations lies, basically, in the failure of Maximum Principle and all the comparison type results based on it. Moreover, it is well known that for higher order problems the lack of the Maximum Principle strongly depends on the kind of boundary conditions imposed on the solutions. In the case of Navier boundary conditions the positivity property is not sensitive to the geometric or topological characteristics of the domain. On the other hand, if we consider Dirichlet boundary conditions, this property is guaranteed only in some special domains, for example in balls taking into account the classical \emph{Boggio's formula} \cite{bog}. Recent results show that it does not hold for certain ellipses \cite{gara,shap} and for squares \cite{coff} but it holds for domains close (in a suitable sense) to the planar disk \cite{gru}, or some non-convex domains \cite{acq}. We refer to the research monograph \cite{gaz} for further details and historical information.\\

The general behavior of arbitrary families of blowing-up solutions to problem \eqref{E1} when $\inf_{\Omega}V>0$ has become understood after (\cite{marpe}, Theorem 1) in the case of Dirichlet boundary condition and after (\cite{lw}, Theorem 1.2) in the case of Navier boundary conditions. It is known that if $u_{\rho}$ is an unbounded family of solutions for which $\rho^{2m}\int_{\Omega}V(x)e^{u_{\rho}}dx$ remains uniformly bounded, then necessarily
\begin{equation}\label{E4}
\lim_{\rho\rightarrow0}\rho^{2m}\int_{\Omega}V(x)e^{u_{\rho}}dx=k\Lambda_{2m}
\end{equation}
for some integer $k\geq1$, where $\Lambda_{2m}:=2^{2m}m!(m-1)!\omega_{2m}$ is a normalization constant and $\omega_{2m}$ denotes the area of the $(2m-1)$-dimensional unit sphere in $\mathbb{R}^{2m}$. Notice that the constant $\Lambda_{2m}$ has a geometric meaning: it is the total $Q$-curvature of the round $2m$-dimensional sphere, see also \cite{marpe,mar} for further information.  Moreover, there are $k$-tuples of distinct points of $\Omega$, the so-called \emph{concentration points}, $(\xi_{1},...,\xi_{k})$, separated at uniformly positive distance from each other and from the boundary $\partial\Omega$ as $\rho\downarrow0$ for which $u_{\rho}$ remains uniformly bounded on $\Omega\setminus\bigcup_{j=1}^{k}B_{\delta}(\xi_{j})$ and $\sup_{B_{\delta}(\xi_{j})}u_{\rho}\rightarrow+\infty$ for any $\delta>0$.
\\An obvious question is the reciprocal, namely construction of solutions to Problem \eqref{E1} with the property \eqref{E4}.
We shall see, firstly, that the location of the concentration points is related to the set of critical points of the finite dimensional functional $\varphi_{k}$ defined explicitly in terms of the potential $V$ and the Green's function $G$ of the  polyharmonic operator with respect to the Dirichlet or Navier boundary condition respectively by
\begin{equation}\label{E17}
\varphi_{k}(\xi)=:-\sum_{j=1}^{k}\big[2\log V(\xi_{j})+H(\xi_{j},\xi_{j})\big]-\sum_{i\neq j}G(\xi_{i},\xi_{j}),
\end{equation}
defined for points $\xi=(\xi_{1},...,\xi_{k})\in\Omega^{k}\setminus\Delta$, where $\Omega^{k}$ denote the cartesian product of $k$ copies of $\Omega$ and $\Delta:=\{\xi_{j}\in\Omega^{k}:\ \xi_{i}=\xi_{j},\ i\neq j\}$ denotes the diagonal in $\Omega^{k}$. For any $\xi\in\overline{\Omega}$, let G=G($x,\xi$) denote the \emph{Green's function} of $(-\Delta)^{m}$ under the appropriate boundary conditions as the unique solution to
\begin{equation}\label{E19}
\left\{ \begin{array}{ll}
          (-\Delta)^{m}G(\cdot,\xi)=\Lambda_{2m}\delta_{\xi}(x)& \mbox{in}\ \Omega\\
          B_{j}G(\cdot,\xi)=0, |j|\leq m-1& \mbox{on}\ \partial\Omega,\
        \end{array} \right.
\vspace{0,2cm}\end{equation}
where $\delta_{\xi}(x):=\delta(x-\xi)$ is the Dirac measure centered at the pole $\xi$. Clearly, the concept of Green's function generalizes that of a fundamental solution. If the differential operator involved has constant coefficients, as in our case, it is in fact often advantageous to think of the Green's function as a perturbation of the fundamental solution. Namely, we decompose $G$ in a singular and regular part as
\begin{equation}\label{E39}
G(x,\xi)=H(x,\xi)+K(x,\xi)
\end{equation}
where $K(x,\xi)=4m\log\frac{1}{|x-\xi|}$ is the fundamental solution of $(-\Delta)^{m}$ in $\mathbb{R}^{2m}$ and the \emph{regular part} of $G$ is a smooth auxiliary polyharmonic function $H(x,\xi)$ so that
\begin{equation}\label{E18}
\left\{ \begin{array}{ll}
(-\Delta)^{m}H(\cdot,\xi) =0& \mbox{in}\ \Omega\\
B_{j}H(\cdot,\xi)=-B_{j}K(\cdot,\xi)& \mbox{on}\ \partial\Omega.\
        \end{array} \right.
\vspace{0,2cm}\end{equation}
Finally, let $H(\xi):=H(\xi,\xi)$ denote the diagonal of the regular part $H$, usually referred as the $\emph{Robin's function}$ of $\Omega$ at $\xi$ and, notice that $H(\xi,\xi)\to -\infty$ as $\xi \to \partial\Omega.$
\\We observe that arbitrary critical points of $\varphi_{k}$ are not all candidates to be concentration points. Since the construction relies on perturbation technique those critical points must also satisfy a sort of nondegeneracy condition.
In \cite{barak2}, the authors solve the question in the case $m=2$, $V\equiv1$ and with Navier boundary conditions. They established that for any \emph{nondegenerate} critical point of $\varphi_{k}$, a family of solutions concentrating at this point as $\rho\downarrow0$ does exist. But, as remarked in \cite{barak2}, their construction, based on a very precise approximation of the actual solution and an application of Banach fixed point theorem, uses nondegeneracy in essential way. This assumption, however, is hard to check in practice and in general not true, an annulus being an obvious example. Next, in \cite{clapp} the authors present a construction of a blowing-up families of solutions under a weaker nondegeneracy assumption of \cite{barak2}, namely that $\varphi_{k}$ has a \emph{topologically nontrivial critical value} and including also a general smooth potential $V$.\\

In the same spirit of \cite{clapp}, we provide sufficient conditions for the existence of multipeak solutions to (\ref{E1}), generalizing all the previous results to the polyharmonic case by also including the Dirichlet boundary condition. We emphasize that all the results presented so far in literature based on the finite dimensional reduction technique are all related to the Navier boundary conditions, see for instance \cite{alar, el, ham, maa, maa2}.
The critical points we will deal with are those that can be captured in a general way with a local min-max characterization. We consider the role of \emph{non-trivial critical values} in an appropriate subset of a functional $\varphi$ in existence of blowing-up solutions to \eqref{E1} in order to relax the stronger and hard-to-check nondegeneracy condition. This local notion of non-trivial critical value was firstly introduced in \cite{delpfel} in the analysis of concentration phenomena in nonlinear Schr\"{o}dinger equation and it is widely used since. More precisely, we consider the following setting. A first observation we make is that in any compact subset of $\widetilde{\Omega}^{k}$, we may define, without ambiguity,
\begin{displaymath}
\varphi(\xi_{1},...,\xi_{k})=-\infty\ \mbox{if}\ \xi_{i}=\xi_{j}\ \mbox{for some}\ i\neq j.
\end{displaymath}
We shall assume that there exists an open subset $\mathcal{D}$ of $\Omega^{k}$ with smooth boundary compactly contained in $\widetilde{\Omega}^{k}$ and such that $\inf_{\mathcal{D}}V>0$. Let $\varphi:\mathcal{D}\rightarrow\mathbb{R}$ be a smooth functional. We will say that \emph{$\varphi$ links in $\mathcal{D}$ at critical level $\mathcal{C}$ relative to $B$ and $B_{0}$} if $B$ and $B_{0}$ are closed subsets of $\overline{\mathcal{D}}$ with $B$ connected and $B_{0}\subset B$ such that the following condition holds.\\
Let us set $\Upsilon$ to be the class of all maps $\Phi \in C(B,\mathcal{D})$ such that there exists an homotopy $\Psi \in C([0,1]\times B,\mathcal{D})$ satisfying:
\begin{displaymath}
\Psi(0,\cdot)=id_{B}, \Psi(1,\cdot)=\Phi, \Psi(t,\cdot)|B_{0}=id_{B_{0}}\ \mbox{for all}\ t \in [0,1].
\end{displaymath}
We assume additionally
\begin{equation}\label{E5}
\sup_{y \in B_{0}}\varphi(y)<\mathcal{C}\equiv \inf_{\Phi \in \Upsilon}\sup_{y \in B}\varphi(\Phi(y)),
\end{equation}
and for every $y \in \partial\mathcal{D}$ such that $\varphi(y)=\mathcal{C}$, there exists a direction vector $\tau_{y}$ tangent to $\partial\mathcal{D}$ at $y$ such that a \emph{transversality condition} holds:
\begin{equation}
\partial_{\tau_{y}}\varphi(y)\equiv\nabla\varphi(y)\cdot\tau_{y}\neq0,
\end{equation}
where $\partial_{\tau_{y}}$ denotes tangential derivative. We observe that this condition can also be stated in the level set framework as $\{y \in \partial\mathcal{D}: \varphi(y)=\mathcal{C}\}=\emptyset$.
Furthermore, we call the min-max value $\mathcal{C}$ a \emph{non-trivial critical level} or a \emph{stable critical level} of the functional $\varphi$ in $\mathcal{D}$. We observe that in the standard language of calculus of variations, the sets $B$, $B_{0}$ \emph{link} in $\mathcal{D}$.
Note that under these assumptions the min-max value $\mathcal{C}$ is a critical value (inside $\Omega$) for $\varphi$ in $\mathcal{D}$, which is in some sense \emph{topologically non-trivial}. Thus these conditions ensure the existence of a critical point of $\varphi$ at level $\mathcal{C}$, i.e. it exists $\overline{y} \in \mathcal{D}$ critical point of $\varphi$ at level $\mathcal{C}$, that is with $\varphi(\overline{y})=\mathcal{C}$. In fact, the condition \eqref{E5} is necessary in order to "seal" $\mathcal{D}$ at level $\mathcal{C}$ so that, exploiting the local linking structure described, by standard deformation argument involving the negative gradient flow of $\varphi$ we are able to provide the presence of such critical point at level $\mathcal{C}$ in $\mathcal{D}$, possibly admitting fully degeneracy. Not only this, any function "$C^{1}$-close" to $\varphi$ (i.e. a $C^{1}$ small perturbation of $\varphi$) inherits such a critical point. Condition \eqref{E5} is a general way of describing a change of topology in the level sets $\{\varphi\leq c\}$ in $\mathcal{D}$ taking place at $c=\mathcal{C}$. The min-max value $\mathcal{C}$ is a critical value if $\mathcal{D}$ is invariant under the negative gradient flow of the functional, while if this is not the case, one can use condition (\ref{E5}) to prevent intersection of the level set $\mathcal{C}$ with the boundary in order to modify the gradient vector field of $\varphi$ near the boundary of $\mathcal{D}$ at the level $\mathcal{C}$ and thus obtain a new vector field with the same stationary points and such that $\mathcal{D}$ is invariant and the functional $\varphi$ is a Lyapunov function for the associated negative flow near the level $\mathcal{C}$. As an example, taking $B=B_{0}=\partial\mathcal{D}$, it is not hard to check that the above conditions hold if
\begin{displaymath}
\inf_{x \in \mathcal{D}}\varphi(x)<\inf_{x \in \partial\mathcal{D}}\varphi(x),\ \mbox{or}\ \sup_{x \in \mathcal{D}}\varphi(x)>\sup_{x \in \partial\mathcal{D}}\varphi(x),
\end{displaymath}
namely the case of (possibly degenerate) general local minimum, maximum points or saddle-points of $\varphi$. The level $\mathcal{C}$ may be taken in these cases respectively as that of the minimum and the maximum of $\varphi$ in $\mathcal{D}$. These holds also if $\varphi$ is "$C^{1}$-close" to a function with a non-degenerate critical point in $\mathcal{D}$. We will show that, for every $k\leq1$, the set $\mathcal{D}:=\{x \in\Omega: dist(x,\partial\Omega)>\delta\}$ has the property \eqref{E5} at a given $\mathcal{C}$, for $\delta$ small enough. This allows us to prove our main result.\newpage
\begin{thm}\label{T1}
\end{thm}
\emph{Let $k\geq1$ and assume that there exists an open subset $\mathcal{D}$ of $\Omega^{k}$ with smooth boundary, compactly contained in $\Omega^{k}$, with $\inf_{\mathcal{D}}V>0$ where $\varphi_{k}$ has a non-trivial critical level $\mathcal{C}$. Then, for $\rho$ small enough, there exists a solution $u_{\rho}$ to \eqref{E1} with $$\lim_{\rho\rightarrow0}\rho^{2m}\int_{\Omega}V(x)e^{u_{\rho}}=k\Lambda_{2m}.$$ Moreover, there is an $k$-tuple $(\xi_{1},...,\xi_{k})\in\mathcal{D}$, such that as $\rho\rightarrow0$ $$\nabla\varphi_{k}(\xi_{1},...,\xi_{k})\rightarrow0,\ \varphi_{k}(\xi_{1},...,\xi_{k})\rightarrow\mathcal{C},$$ for which $u_{\rho}$ remains uniformly bounded on $\Omega\setminus\bigcup_{j=1}^{k}B_{\delta}(\xi_{j})$ and, for any $\delta>0$, $\sup_{B_{\delta}(\xi_{j})}u_{\rho}\rightarrow+\infty$.}\\
\\Notice that for the two-dimensional version of problem \eqref{E1} many authors, using different perturbation techniques, have constructed solutions which admit a concentration behavior with a prescribed \emph{bubble} profile as indicated by the blow-up analysis in the seminal papers \cite{12,31,32,38}. In \cite{barak}, providing that $\Omega$ is not simply connected, the authors showed that for any non-degenerate critical point of the reduced function a sequence $u_{\rho}$ of solutions is constructed, which converges to a function $u^{\ast}=8\pi\sum_{i=1}^{m}G(x,\xi_{i})$ in $C^{2,\alpha}_{loc}(\overline{\Omega}\setminus S)$.
After, in \cite{delpi,esp} the authors, independently, generalize this result by relaxing the assumption of non degeneracy into a weaker stability assumption, the \emph{topologically\ nontrivial\ critical\ value} used also in our paper.\\
\\The paper is organized as follows. Section \ref{w} is devoted to describing a first approximation for the solution and to estimating the error. The predicted solutions are found as a small perturbation of this initial approximation. Furthermore, problem \eqref{E1} is written as a fixed point problem, involving a linear operator. In Section \ref{an} we study the bounded invertibility of the linearized operator in a suitable $L^{\infty}$-weighted space. In Section \ref{pro} we solve a projected non-linear problem. In this Section, in order to include Dirichlet boundary conditions we are able to avoid the use of Maximum Principle in the proof of Lemma \ref{aa}. This result represents the major technical difference with respect to \cite{clapp}. In Sections \ref{S1} and \ref{S10} we show that solving the entire non-linear problem reduces to finding critical points of a certain functional and this correspond to adjust variationally the location and the high of the bubbles. Finally, Section \ref{S20} is devoted to the proof of the main result.\\

\section{A first approximation of solution}\label{w}
In this section we construct a reasonably good approximation $U$ for a solution of \eqref{E1}, we give an expansion of its associated energy and then we estimate the error of such approximation in appropriate norms. Our construction relies in an essential way in two main steps:
\begin{itemize}
\item on the choice of a suitable family of approximating solution,
\item on the invertibility, in some sense, of the linearized operator evaluated at such approximating solution.
\end{itemize}
A useful tool involved in this construction is, in fact, concerned with the classification of entire solutions to the higher order Liouville equation which represents, essentially, the associate limit problem to \eqref{E1}.
Given a constant $Q\in\mathbb{R}$, we consider the following \emph{limit profile} problem
\begin{equation}
\left\{ \begin{array}{ll}
          (-\Delta)^{m} u=\ Qe^{2mu}& \mbox{in}\ \mathbb{R}^{2m}, m>1\\
          e^{2mu}\in L^{1}(\mathbb{R}^{2m}),&\\
          u(x)=o(|x|^{2})& \mbox{as}\ |x|\rightarrow\infty.
        \end{array} \right.
\vspace{0,2cm}\end{equation}
In \cite{mar} it has been proved that all the solutions of \eqref{E1} are radially symmetric about some point $\xi\in\mathbb{R}^{2m}$ and are all given by the one-parameter family of functions
\begin{equation}\label{E6}
U_{\delta,\xi}(x)= \log\frac{\alpha_{2m}Q\delta^{2m}}{(\delta^{2}+|x-\xi|^{2})^{2m}}
\end{equation}
for any free parameter $\delta>0$, $\xi\in\mathbb{R}^{2m}$ and where $\alpha_{2m}:=2^{2m+1}m(m-1)!$ if and only if the constant $Q$ is strictly positive as pointed out in \cite{mar2}, if $Q\leq0$ and $m>1$ there are no standard solutions and there are no solutions presenting a nice behavior at infinity. We call the functions of this form $\emph{standard}$ (or $\emph{geometric}$) solutions. Notice that given a solution $u$ to $(-\Delta)^{m} u=\ Qe^{2mu}\ \mbox{in}\ \mathbb{R}^{2m}, m>1$ and $\lambda>0$, the new function defined as $v:=u-\frac{1}{2m}\log\lambda$ solves $(-\Delta)^{m} v=\ \lambda Qe^{2mv}\ \mbox{in}\ \mathbb{R}^{2m}$, hence what really matters is just the sign of $Q$. There is no loss of generality if one assumes that $Q\in\{0,\pm(2m-1)!\}$. On the other hand, the particular positive choice $Q=(2m-1)!$ has the advantage of being exactly the $Q$-curvature of the round sphere $S^{2m}$. This fact implies that the standard solutions has the geometric property that $e^{2u}g_{\mathbb{R}^{2m}}=(\pi^{-1})^{\star}g_{S^{2m}}$, where $\pi:S^{2m}\rightarrow\mathbb{R}^{2m}$ is the stereographic projection, then by M\"{o}bius transformations actually give us the large family of solutions \eqref{E6}. We recall that in dimension 2, by employing the method of moving plane, in \cite{15} the authors were able to classify $\emph{all}$ the radial solutions of
\begin{equation}
\left\{ \begin{array}{ll}
          (-\Delta)u=\ e^{2u}\ \mbox{in}\ \mathbb{R}^{2}, \\
          e^{2u}\in L^{1}(\mathbb{R}^{2})\\
        \end{array} \right.
\vspace{0,2cm}\end{equation}
as standard solutions without any other condition. Conversely, if $m>1$ a careful study of radial solutions shows that there are solutions which do not come from the smooth function on $S^{2m}$ through the stereographic projection. Thus, to overcome this difficulty, first in \cite{wxu} the author added a constraint in the behavior at infinity of the solution in order to have the same classification obtained in $m=1$ case. Notice that their proof appears to be overly simplified, for instance their Lemma 2.2 is not conclusive. In (\cite{mar}, Theorem 2) the author give a proof of this result with all the details and develop some criteria to characterize of the non-standard solutions. Of particular importance is the following criterium: If a solution $u$ is non-standard, then there exist $1\leq j\leq m-1$  and a constant $a\neq0$ such that $$\lim_{|x|\rightarrow\infty}\Delta^{j}u(x)=a.$$ Notice that a similar property can be shown for \emph{every} solution in the negative case.
Moreover, the standard solutions satisfy some nice properties:
\begin{enumerate}
  \item $\int_{\mathbb{R}^{2m}} e^{U_{\delta,\xi}(x)}dx=\Lambda_{2m}$ (\emph{mass quantization}),
  \item $e^{U_{\delta,\xi}(x)}\rightharpoonup\Lambda_{2m}\delta_{\xi}$ in measure sense as $\delta\rightarrow0$ ($\emph{concentration property}$),
  \item given a small $\delta_{0}>0$, $\sup_{B_{\delta_{0}}(\xi)}U_{\delta,\xi}(x)\rightarrow+\infty$ as $\delta\rightarrow0$ and $U_{\delta,\xi}$ remains uniformly bounded on $\Omega\setminus{B_{\delta_{0}}(\xi)}$ \ ($U_{\delta,\xi}(x)\ \emph{blows up}\ at\ \xi$ and remains bounded away from it),
  \item $U_{\delta,\xi}(x)\rightarrow-\infty$ as $|x|\rightarrow+\infty.\\$
\end{enumerate}

Due to these all properties we shall use $U_{\delta,\xi}$, suitable scaled and projected, as fundamental $\emph{building blocks}$ to construct an approximate solution of \eqref{E1} around $\xi$.
Given $k$ a positive integer, let us consider $k$ distinct and well-separated points $\xi_{1},...,\xi_{k} \in \Omega$ where the spikes are meant to take place with $V(\xi_{i})>0$ for all $i=1,...,k$ because, as we said before, for $Q\leq0$ and $m>1$ there are no standard solutions at all. As we will see, a convenient set to select these points is
\begin{equation}
\mathcal{M}:=\big\{\xi\in\Omega^{k}: \mbox{dist}(\xi_{i},\partial\Omega)\geq2\delta_{0}\ \forall i=1,...,k, \min_{i\neq j}|\xi_{i}-\xi_{j}|\geq2\delta_{0}\big\}
\end{equation}
where $\delta_{0}>0$ is a small but fixed number. From now on we fix $\xi=(\xi_{1},...,\xi_{k})\in \mathcal{M}$.
Such a definition is motivated by a blow-up analysis performed in.
Given $\varepsilon>0$, we define for $x\in\mathbb{R}^{2m}$ around each $\xi_{i}$
\begin{equation}
u_{\varepsilon}(x):= 2m \log\frac{(1+\varepsilon^{2})}{(\varepsilon^{2}+|x-\xi_{i}|^{2})V(\xi_{i})}
\end{equation}
so that $u_{\varepsilon}$ is a rotationally symmetric solution to
\begin{equation}\label{E7}
(-\Delta)^{m} u=\rho^{2m} V(\xi_{i}) e^{u}\ \mbox{in}\ \mathbb{R}^{2m}
\end{equation}
with
\begin{equation}\label{E16}
\rho^{2m}= \frac{\alpha_{2m}(2m-1)!\varepsilon^{2m}}{(1+\varepsilon^{2})^{2m}},
\end{equation}
that is, $\rho\sim\varepsilon$ as $\varepsilon\rightarrow0$.
Now, let us notice that equation \eqref{E7} is scale invariant under some dilation in the following sense: if $u_{\varepsilon}$ is a solution of \eqref{E7} and $\mu_{i}>0$, $i=1,...,k$, are $k$ positive parameters to be selected properly later on, then $u_{\varepsilon}(\mu_{i}\cdot)+2m\log\mu_{i}$ ia also a solution of \eqref{E7}. Thus, with this observation in mind, we define for all $\mu_{i}>0$
\begin{equation}\label{E12}
u_{i}(x):= 2m\log\frac{\mu_{i}(1+\varepsilon^{2})}{(\mu_{i}\varepsilon^{2}+|x-\xi_{i}|^{2})V(\xi_{i})}
\end{equation}
so that the solution $u_{\varepsilon}$ will look, near each $\xi_{i}$, like $u_{i}$ for certain $\varepsilon$-independent parameters $\mu_{i}$.
\begin{oss}
\end{oss}
\begin{enumerate}
  \item Note also that $u_{i}$ uniquely solves \eqref{E7} in entire $\mathbb{R}^{2m}$ with $\rho\sim\varepsilon$ as $\varepsilon\rightarrow0$, in fact, up to scaling and translation invariance, these are the only solutions with finite energy condition and a nice behavior at infinity.
  \item A considerable characteristic of elliptic PDEs with nonlinearity of exponential type in the $\emph{right}$ dimension (i.e. $(-\Delta)^{m}\ \mbox{in}\ \mathbb{R}^{2m}$) is that all the bubble solutions have proportional heights with coefficient of proportionality $\mu_{i}$.
\end{enumerate}
We would like to take $\sum_{i=1}^{k}u_{i}$ as first approximation to a solution of \eqref{E1}.
Since, obviously, $B_{j}u_{i}$, where $|j|\leq m-1$ and $i=1,...,k$, are not zero on the boundary $\partial\Omega$, we perform a suitable polyharmonic correction to modify $\sum_{i=1}^{k}u_{i}$ in order to fit in the boundary conditions.
Let $H_{i}(x)$ be the smooth unique solution of
\begin{equation}
\left\{\begin{array}{ll}
          (-\Delta)^{m}H_{i}(x) =0& \mbox{in}\ \Omega\\
          B_{j}H_{i}(x)=B_{j}u_{i}(x), |j|\leq m-1& \mbox{on}\ \partial\Omega.\
        \end{array}\right.
\vspace{0,2cm}\end{equation}
We define our first approximation $U(\xi)$ as
\begin{equation}\label{E14}
U(\xi):=\sum_{i=1}^{k}U_{i},\ U_{i}:=u_{i}+H_{i}
\end{equation}
which now, by definition, satisfies the boundary conditions.
We want to find out the asymptotic behavior of $U_{i}$ away from $\xi_{i}$ and around $\xi_{i}$.\\ First, we recall the standard regularity statement and an \emph{a priori} estimate for higher order linear problems in the $L^{p}$ space framework developed by Agmon, Douglis and Nirenberg \cite{adn} about 50 years ago. Before this, we remark that, in general, the choice of the linear boundary operator $B_{j}$'s for polyharmonic problems is not completely free. In order to obtain a priori estimates and, in turn, existence and uniqueness results we need to impose a certain algebraic constraint, the so-called \emph{Shapiro-Lopatinksi$\check{i}$ complementary condition}. For any $j$, let $B_{j}'$ denote the highest order part of $B_{j}$ which is precisely of order $m_{j}$. For every point $x\in\Omega$, let $\nu(x)$ denote the normal unit vector. We say that the complementing condition holds for $B_{j}u=g_{j}$ on $\partial\Omega$ if, for any nontrivial tangential vector $\tau(x)$, the polynomials in $t$ $B_{j}'(x;\tau+t\nu)$ are linearly independent modulo the polynomial $(t-i|\tau|)^{m}$ with $i$ the unit imaginary number. We remark that both Dirichlet and Navier boundary conditions satisfy this condition. For example, the homogeneous Neumann boundary conditions $\Delta u=\partial_{\nu}(\Delta u)=0$ on $\partial\Omega$ for the biharmonic equation $\Delta^{2}u=0$ in $\Omega$ not satisfy it.
\begin{lm}\label{L1}
\end{lm}
\emph{Let $\Omega$ a bounded domain in $\mathbb{R}^{2m}$. Let $1<p<\infty$ and fix an integer $q\geq 2m$. Assume $\partial\Omega\in C^{q}$. Consider the following problem
\begin{equation}\label{E8}
\left\{\begin{array}{ll}
          (-\Delta)^{m}u(x) =f& \mbox{in}\ \Omega\\
          B_{j}u(x)=g_{j}(x), |j|\leq m-1& \mbox{on}\ \partial\Omega,\
        \end{array}\right.
\vspace{0,2cm}\end{equation}
for all $f\in W^{q-2m,p}(\Omega)$, for all $g_{j}\in W^{q-m_{j}-\frac{1}{p},p}(\partial\Omega)$ with $m_{j}\in\mathbb{N}$ the maximal order of derivatives of $B_{j}$ and assume that $B_{j}$ satisfy the Shapiro-Lopatinksi$\check{i}$ condition. Then \eqref{E8} admits a unique strong solution $u\in W^{q,p}(\Omega)$. Moreover, there exist a constant $c=c(|\Omega|,k,m,B_{j})>0$ independent of $f$ and of the $g_{j}$'s, such that the following a priori estimate holds
\begin{equation}
\|u\|_{W^{q,p}(\Omega)}\leq c\bigg(\|f\|_{W^{q-2m,p}(\Omega)}+\sum_{j=1}^{m}\|g_{j}\|_{W^{q-m_{j}-\frac{1}{p},p}(\partial\Omega)}\bigg).
\end{equation}}
Now, arguing as in the classical \cite{rey}, we obtain the following crucial characterization in which we asymptotically expands $U_{i}$ in $\Omega$.
\begin{lm}\label{a}
\end{lm}
\emph{Assume $\xi\in\mathcal{M}$. Then we have, as $\varepsilon\rightarrow0$ :
\begin{enumerate}
\item uniformly in $\bar{{\Omega}}$ in the $C^{2m-1,\alpha}_{loc}$-sense
\begin{equation}\label{E9}
H_{i}(x)=H(x,\xi_{i})-2m\log\mu_{i}(1+\varepsilon^{2})+\log V(\xi_{i})+O(\mu_{i}^{2}\varepsilon^{2}),
\end{equation}
\item uniformly in the region $|x-\xi_{i}|\geq\delta_{0}$
\begin{equation}\label{E10}
 u_{i}(x)=2m\log \mu_{i}(1+\varepsilon^{2})-\log V(\xi_{i})-4m\log|x-\xi_{i}|+O(\mu_{i}^{2}\varepsilon^{2}).
\end{equation}
\item In particular, in this region, as we expect by the asymptotic analysis,
\begin{equation}\label{E11}
U_{i}(x)=G(x,\xi_{i})+O(\mu_{i}^{2}\varepsilon^{2})
\end{equation}
\end{enumerate}
where $\alpha$ is an arbitrary H\"{o}lder exponent, $G(x,\xi)$ is the Green's function and $H(x,\xi)$ its regular part.
}\\$\underline{Proof}$\\
Let us prove \eqref{E9}. Define $z(x):=H_{i}(x)-H(x,\xi_{i})+2m\log\mu_{i}(1+\varepsilon^{2})-\log V(\xi_{i})$, then $z\in\mathcal{H}_{m}(\Omega)$ and it satisfies
\begin{displaymath}
\left\{\begin{array}{ll}
          (-\Delta)^{m}z(x)=0& \mbox{in}\ \Omega\\
          z(x)=-u_{i}(x)-4m\log|\cdot-\xi_{i}|+2m\log\mu_{i}(1+\varepsilon^{2})-\log V(\xi_{i})& \mbox{on}\ \partial\Omega\\
          B_{j}z(x)=B_{j}u_{i}(x), |j|\leq m-1& \mbox{on}\ \partial\Omega.\
        \end{array}\right.
\vspace{0,2cm}\end{displaymath}
Since $B_{j}z(x)=O(\mu_{i}^{2}\varepsilon^{2})$ with $|j|\leq m-1$, using the standard $L^{p}$-estimates in Lemma \ref{L1} and the Sobolev embedding theorem we get $z(x)=O(\mu_{i}^{2}\varepsilon^{2})$ uniformly in $\overline{\Omega}$ in the $C^{2m-1,\alpha}_{loc}$-sense.\\
The point \eqref{E10} follows, easily, from the definition of $u_{i}$ in \eqref{E12}. Finally, combining \eqref{E9} and \eqref{E10} we get \eqref{E11}.
\begin{flushright}
$\Box$
\end{flushright}
While $u_{i}(x)$ is a good approximation to a solution of \eqref{E1} near $\xi_{i}$, $U(\xi)$ is a good approximation far from these points but, unfortunately, it is not good enough for our constructions close to these points. This is the reason why we need to further adjust our ansatz. To do this, we need to improve that the remainder term  $U-u_{i}(x)=(H_{i}+\sum_{k\neq i}u_{k})\sim0$, that is it vanishes at main order near $\xi_{i}$ as $\varepsilon\rightarrow0$. Using Lemma $\ref{a}$, we can achieve this through the following precise selection of the parameters $\mu_{i}$
\begin{equation}\label{E15}
2m\log\ \mu_{i}=H(\xi_{i},\xi_{i})+\log V(\xi_{i})+\sum_{i\neq j}G(\xi_{j},\xi_{i}).
\end{equation}
We thus fix $\mu_{i}\ \emph{a priori}$ as a function of $\xi\in\mathcal{M}$, $\mu_{i}=\mu_{i}(\xi)\ \forall i=1,...,k$.
Since $\xi\in\mathcal{M}$, for some constant $c>0$,
\begin{equation}
\frac{1}{c}\leq\mu_{i}\leq c\ \mbox{for\ all}\ i=1,...,k.
\end{equation}
and, thus, instead of $O(\mu_{i}^{2}\varepsilon^{2})$ we can write $O(\varepsilon^{2})$ in all the following asymptotic expansions.\\
Now, let us denote with
\begin{equation}
\Omega_{\varepsilon}:=\varepsilon^{-1}\Omega=\{y\in \mathbb{R}^{2m}:\varepsilon y \in \Omega\}
\end{equation}
the $\emph{expanded domain}$ with the scaling measure property $|\Omega_{\varepsilon}|=\varepsilon^{-2m}|\Omega|$ and
\begin{equation}
\xi_{i}':=\varepsilon^{-1}\xi_{i},
\end{equation}
the $\emph{expanded variables}$.
A useful observation is that $u$ satisfies \eqref{E1} if and only if $w(y):=u(\varepsilon y)+ 2m\log \rho\varepsilon$ solves the equivalent problem
\begin{equation}\label{E13}
\left\{ \begin{array}{ll}
          (-\Delta)^{m} w= V(\varepsilon y) e^{w}& \mbox{in}\ \Omega_{\varepsilon}\\
          w=2m\log \rho\varepsilon& \mbox{on}\ \partial\Omega_{\varepsilon}\\
          B_{j}w=0, |j|\leq m-1& \mbox{on}\ \partial\Omega_{\varepsilon}.\
        \end{array} \right.
\vspace{0,2cm}\end{equation}
Let us define $W(y):=U(\varepsilon y)+2m\log\rho\varepsilon$ as the first approximation of \eqref{E13} with $U$ our approximation solution \eqref{E14}.
We want see how $W$ behaves, namely, we want to measure how well it solves \eqref{E13}. It is convenient to do so in terms of the following $L^{\infty}$-weighted norm defined $\mbox{for all}\ w \in L^{\infty}(\Omega_{\varepsilon})$
\begin{equation}\label{Ea}
\|w\|_{\ast}:=\sup_{y \in \Omega_{\varepsilon}}\bigg(\sum_{i=1}^{k}\frac{1}{(1+|y-\xi_{i}'|)^{4m-1}}+\varepsilon^{2m}\bigg)^{-1}|w(y)|.
\end{equation}
Thus, as anticipated, in the next lemma we measure the size of the error of approximation $R(y):=(-\Delta)^{m} W-V(\varepsilon y) e^{W}$, proving that although is not zero it is small in the sense of the norm defined above. We remark that the suitable choice of the parameters $\mu_{i}$ in \eqref{Ea} is done just to make the error term small.
\begin{lm}\label{b}
\end{lm}
\emph{Assume $\xi\in\mathcal{M}$. Then we have for any $y\in\Omega_{\varepsilon}$, as $\varepsilon\rightarrow0$:
\begin{equation}
\|R\|_{\ast}\leq C\varepsilon
\end{equation}
where $C$ denotes a generic positive constant independent of $\varepsilon$ and $\xi$.
}\\$\underline{Proof}$\\
Let us assume first $|y-\xi_{n}'|<\frac{\delta_{0}}{\varepsilon}$, for some index $n$. We have
\begin{displaymath}
(-\Delta)^{m}W(y)=(-\Delta)^{m}(U(\varepsilon y)+2m\log\rho\varepsilon)=(-\Delta)^{m}U(\varepsilon y)=\varepsilon^{2m}[(-\Delta)^{m}U](\varepsilon y)
\end{displaymath}
because the higher order Laplace operator is linear and it scales like a homogeneous polynomial of degree $2m$. Now, recalling \eqref{E14} and the fact that $H_{i}\in\mathcal{H}_{m}(\Omega_{\varepsilon})$ we get
$$\varepsilon^{2m}[(-\Delta)^{m}U](\varepsilon y)=\varepsilon^{2m}[(-\Delta)^{m}u_{i}](\varepsilon y)=\varepsilon^{2m}\rho^{2m}\sum_{i}V(\xi_{i})e^{u_{i}(\varepsilon y)}$$ $$=\varepsilon^{2m}\rho^{2m}V(\xi_{n})\frac{\mu_{n}^{2m}(1+\varepsilon^{2})^{2m}}{(\mu_{n}^{2}\varepsilon^{2}+|\varepsilon y-\xi_{n}|^{2})^{2m}V(\xi_{n})}+\varepsilon^{2m}\rho^{2m}\sum_{i\neq n}V(\xi_{i})e^{u_{i}(\varepsilon y)}$$ $$=\frac{\varepsilon^{4m}\alpha_{2m}(2m-1)!}{(1+\varepsilon^{2})^{2m}}
\frac{\mu_{n}^{2m}(1+\varepsilon^{2})^{2m}}{\varepsilon^{4m}(\mu_{n}^{2}\varepsilon^{2}+|y-\xi_{n'}|^{2})^{2m}V(\xi_{n})}
+\varepsilon^{2m}\rho^{2m}\sum_{i\neq n}V(\xi_{i})e^{u_{i}(\varepsilon)}
$$
\begin{displaymath}
=\alpha_{2m}(2m-1)!\frac{\mu_{n}^{2m}(1+\varepsilon^{2})^{2m}}{(\mu_{n}^{2}\varepsilon^{2}+|y-\xi_{n}'|^{2})^{2m}V(\xi_{n})}
+\varepsilon^{2m}\rho^{2m}\sum_{i\neq n}V(\xi_{i})e^{u_{i}(\varepsilon y)}
\end{displaymath}
\begin{equation}
=\alpha_{2m}(2m-1)!\frac{\mu_{n}^{2m}(1+\varepsilon^{2})^{2m}}{(\mu_{n}^{2}\varepsilon^{2}+|y-\xi_{n}'|^{2})^{2m}V(\xi_{n})}
+O(\varepsilon^{4m})
\end{equation}
because of $i\neq n$ means that $y$ is away from the points $\xi_{i}'$ and then we should use \eqref{E10}.
Let us estimate $V(\varepsilon y)e^{W(y)}$. By \eqref{E9} we have
\begin{displaymath}
H_{k}(x)=H(\xi_{n},\xi_{n})-2m\log\mu_{n}+\log V(\xi_{n})+O(\varepsilon^{2})+O(|x-\xi_{n}|)
\end{displaymath}
and by the definition of $\mu_{i}$'s \eqref{E15} we have
\begin{displaymath}
=-\sum_{i\neq n}G(\xi_{i},\xi_{n})+O(\varepsilon^{2})+O(|x-\xi_{n}|)
\end{displaymath}
and, if $i\neq n$, by \eqref{E11} we have
\begin{displaymath}
U_{i}(x)=G(\xi_{i},\xi_{i})+O(\varepsilon^{2})+O(|x-\xi_{n}|).
\end{displaymath}
Then, by adding, we have
\begin{equation}\label{E55}
H_{n}(x)+\sum_{i\neq n}U_{i}(x)=O(\varepsilon^{2})+O(|x-\xi_{n}|).
\end{equation}
Therefore,
\begin{displaymath}
V(\varepsilon y)e^{W(y)}=V(\varepsilon y)(\rho\varepsilon)^{2m}e^{U(\varepsilon y)}=V(\varepsilon y)(\rho\varepsilon)^{2m}\exp\bigg(u_{n}(\varepsilon y)+H_{n}(\varepsilon y)+\sum_{i\neq n}U_{i}(\varepsilon y)\bigg)
\end{displaymath}
and using \eqref{E55},
\begin{eqnarray*}
&=&V(\varepsilon y)(\rho\varepsilon)^{2m}\exp\bigg(u_{n}(\varepsilon y)+O(\varepsilon^{2})+O(\varepsilon|y-\xi_{m}')\bigg)\\
&=&V(\varepsilon y)(\rho\varepsilon)^{2m}\frac{\mu_{n}^{2m}(1+\varepsilon^{2})^{2m}}{(\mu_{k}^{2}\varepsilon^{2}+|y-\xi_{n}'|^{2})V(\xi_{n})}
\bigg(1+O(\varepsilon^{2})+O(\varepsilon|y-\xi_{n}'|)\bigg)
\end{eqnarray*}
and, finally,
\begin{displaymath}
=\alpha_{2m}(2m-1)!\frac{\mu_{n}^{2m}}{(\mu_{n}^{2}+|y-\xi_{n}'|^{2})^{2m}}\bigg(1+O(\varepsilon^{2})+O(\varepsilon|y-\xi_{n}'|)\bigg)
\end{displaymath}
because by first order Taylor expansion we have
\begin{displaymath}
\frac{V(\varepsilon y)}{V(\xi_{n})}=1+O(\varepsilon|y-\xi_{n}'|).
\end{displaymath}
We can conclude that in this region
\begin{displaymath}
|R(y)|\leq C\varepsilon\frac{1}{(1+|y-\xi_{n}'|)^{4m-1}}
\end{displaymath}
and thus $\|R(y)\|_{\ast}=O(\varepsilon)$ as $\varepsilon\rightarrow0$ if $y\in B_{\frac{\delta_{0}}{\varepsilon}}(\xi_{n}')$.\\
Now, we assume $|y-\xi_{i}'|\geq\frac{\delta_{0}}{\varepsilon}$ for all $i=1,...,k$, scaling back this is equivalent to say $|x-\xi_{i}|\geq\delta_{0}$ so, using \eqref{E11},
\begin{displaymath}
(-\Delta)^{m}W(y)=\varepsilon^{2m}(-\Delta)^{m}U(\varepsilon y)=\varepsilon^{2m}(-\Delta)^{m}(G(\varepsilon y,\xi_{i})+O(\varepsilon^{2})=O(\varepsilon^{4m})
\end{displaymath}
and $V(\varepsilon y)e^{W(y)}=O(\varepsilon^{4m})$. Hence, in this region,
\begin{displaymath}
R(y)=O(\varepsilon^{4m})
\end{displaymath}
so that, finally,
\begin{displaymath}
\|R(y)\|_{\ast}=O(\varepsilon)\ \mbox{as}\ \varepsilon\rightarrow0\ \mbox{in all}\ \Omega_{\varepsilon}
\end{displaymath}
and the proof is concluded.
\begin{flushright}
$\Box$
\end{flushright}
Observe that problem \eqref{E1} has a $\emph{variational structure}$, in the sense that (weak) solutions for \eqref{E1} correspond to critical points of the following nonlinear energy functional
\begin{equation}\label{E56}
J_{\rho}(u):=\frac{1}{2}\int_{\Omega}|(-\Delta)^{\frac{m}{2}}u|^{2}dx-\rho^{2m}\int_{\Omega}V(x)e^{u}dx
\end{equation}
where $u$ acts on the Hilbertian Sobolev space $H=\{u\in H^{m}(\Omega)\ \mbox{such that}\ B_{j}u=0\ \mbox{on}\ \partial\Omega, |j|\leq m-1\}$, the natural settings for the boundary operator involved.
In particular, in the case of Dirichlet boundary conditions we have
\begin{displaymath}
H= H_{0}^{m}(\Omega),
\end{displaymath}
i.e. $H$ is defined as the completion of $C_{0}^{\infty}(\Omega)$ with respect to the standard Sobolev norm $\|(-\Delta)^{\frac{m}{2}}u\|_{L^{2}(\Omega)}$ defined below, while in the case of Navier boundary conditions we have
\begin{displaymath}
H=H_{N}:=\bigg\{u\in H^{m}(\Omega)\ \mbox{such that}\ (-\Delta)^{j}u=0\ \mbox{on}\ \partial\Omega, j=0,1,...,\bigg[\frac{m-1}{2}\bigg]\bigg\}.
\end{displaymath}
$H$ inherits the Hilbert space structure from $H^{m}(\Omega)$ when endowed with the scalar product
\begin{displaymath}
(u,v)_{H}:=
\left \{ \begin{array}{ll} \displaystyle\int_{\Omega}(-\Delta)^{k}u (-\Delta)^{k}v\ dx& \mbox{if}\ m=2k\\
\displaystyle\int_{\Omega}\nabla(-\Delta)^{k}u\nabla(-\Delta)^{k}v\ dx& \mbox{if}\ m=2k+1
        \end{array}\right.
\vspace{0,2cm}\end{displaymath}
and where we use the following notation:
\begin{displaymath}
(-\Delta)^{\frac{m}{2}}u:=
\left \{ \begin{array}{ll} (-\Delta)^{k}u& \mbox{if}\ m=2k\\
\nabla(-\Delta)^{k}u& \mbox{if}\ m=2k+1.
        \end{array}\right.
\vspace{0,2cm}\end{displaymath}
The scalar product defined above induces the following corresponding norm on $H$:
\begin{displaymath}
\|u\|_{H}:=(u,v)_{H}^{\frac{1}{2}}=
\left \{ \begin{array}{ll} \|(-\Delta)^{k}u\|_{L^{2}(\Omega)}\in \mathbb{R}& \mbox{if}\ m=2k\\
\|\nabla(-\Delta)^{k}u\|_{L^{2}(\Omega)}\in \mathbb{R}^{2m}& \mbox{if}\ m=2k+1.
        \end{array}\right.
\vspace{0,2cm}\end{displaymath}
For general $p\in(1,+\infty)$, one has the choice of taking the $L^{p}$-version of this norm or the equivalent Sobolev norm $\|u\|_{W_{0}^{m,p}(\Omega)}:=\|D^{m}u\|_{L^{p}(\Omega)}$.\\
\\We end this section giving an asymptotic estimate of the $\emph{reduced energy}$ $J_{\rho}(U)$, that is the energy functional defined above calculated along our ansatz $U$. Instead of $\rho$, we use the parameter $\varepsilon$, related to $\rho$ by \eqref{E16}, to obtain the following result:
\begin{lm}\label{c}
\end{lm}
\emph{With the election of $\mu_{i}$'s given by \eqref{E15},
\begin{equation}\label{E57}
J_{\rho}(U)=b_{m}\varphi_{k}(\xi)+4mb_{m}k|\log\varepsilon|-4b_{m}k+O(\varepsilon)
\end{equation}
where $\xi\in\mathcal{M}$, $\varphi_{k}$ is the functional defined in \eqref{E17} and $b_{m}:=\frac{1}{2}\alpha_{2m}(m-1)!\pi^{m}$.
}\\$\underline{Proof}$\\
We consider the energy evaluated at $U$, that is
\begin{displaymath}
J_{\rho}(U)=\frac{1}{2}\sum_{i=1}^{k}\int_{\Omega}|(-\Delta)^{\frac{m}{2}}U_{i}|^{2}+\frac{1}{2}\sum_{i\neq j}\int_{\Omega}|(-\Delta)^{\frac{m}{2}}U_{i}||(-\Delta)^{\frac{m}{2}}U_{j}|-\rho^{2m}\int_{\Omega}V(x)e^{U}
\end{displaymath}
$\equiv I_{1}+I_{2}+I_{3}$. We will first evaluate the quadratic part $I_{1}+I_{2}$.\\
Let $i$ be fixed, we note that $(-\Delta)^{m}U_{i}(x)=(-\Delta)^{m}(u_{i}+H_{i})=(-\Delta)^{m}u_{i}(x)=\rho^{2m}V(\xi_{i})e^{u_{i}}$ in $\Omega$ and $B_{j}U_{i}=0$ on $\partial\Omega$ where $i=1,...,k$ and $j=0,...,m-1$. Then, integrating by parts,
\begin{eqnarray*}
I_{1}&=&\frac{1}{2}\rho^{2m}\sum_{i=1}^{k}V(\xi_{i})\int_{\Omega}e^{u_{i}}U_{i},\\
I_{2}&=&\frac{1}{2}\rho^{2m}\sum_{i\neq j}V(\xi_{i})\int_{\Omega}e^{u_{i}}U_{j}.
\end{eqnarray*}
First we expand $I_{1}$, using \eqref{E9} and \eqref{E16} and recalling that $H_{i}\in\mathcal{H}_{m}(\Omega)$ we have
\begin{eqnarray*}
I_{1}&=&\frac{1}{2}\rho^{2m}\sum_{i=1}^{k}V(\xi_{i})\int_{\Omega}e^{u_{i}}(u_{i}+H_{i})\\
&=&\frac{1}{2}\alpha_{2m}(2m-1)!\sum_{i=1}^{k}\int_{\Omega}
\frac{(\mu_{i}\varepsilon)^{2m}}{(\mu_{i}^{2}\varepsilon^{2}+|x-\xi_{i}|^{2})^{2m}}\bigg\{-2m\log(\mu_{i}^{2}\varepsilon^{2}\\
&+&|x-\xi_{i}|^{2})+H(x,\xi_{i})+O(\varepsilon^{2})\bigg\}.
\end{eqnarray*}
Now let $\delta_{0}>0$ be small and fixed, independent of $\rho$. Let us define the change of variables $x=\xi_{i}+\mu_{i}\varepsilon y$, where $x\in B_{\delta_{0}}(\xi_{i})$ and $y\in \frac{B_{\delta_{0}}(\xi_{i})-\xi_{i}}{\mu_{i}\varepsilon}\equiv B_{\frac{\delta_{0}}{\mu_{i}\varepsilon}}(0)$ with the relative volume element $dx=(\mu_{i}\varepsilon)^{2m} dy$.
As usual, we can split the above integral on $\Omega$ into two pieces, a sum of one integral on $B_{\frac{\delta_{0}}{\mu_{i}\varepsilon}}(0)$ and one on its complement. Since the piece on the complement is small, precisely an $O(\varepsilon^{2m})$, it is enough expand the integral only on $B_{\frac{\delta_{0}}{\mu_{i}\varepsilon}}(0)$.
\begin{eqnarray*}
I_{1}&=&\frac{1}{2}\alpha_{2m}(2m-1)!\sum_{i=1}^{k}\int_{B_{\frac{\delta_{0}}{\mu_{i}\varepsilon}}(0)}
\frac{(\mu_{i}\varepsilon)^{4m}dy}{(\mu_{i}^{2}\varepsilon^{2}+\mu_{i}|\varepsilon y|^{2})^{2m}}\bigg\{-2m\log(\mu_{i}^{2}\varepsilon^{2}+|\mu_{i}\varepsilon y|^{2})\\
&+&H(\xi_{i}+\mu_{i}\varepsilon y,\xi_{i})+O(\varepsilon^{2})\bigg\}+O(\varepsilon^{2m}).
\end{eqnarray*}
Bringing out the term $(\mu_{i}^{2}\varepsilon^{2})^{2m}$ and using first order Taylor expansion:
\begin{displaymath}
H(\xi_{i}+\mu_{i}\varepsilon y,\xi_{i})=H(\xi_{i},\xi_{i})+O(\mu_{i}|\varepsilon y|),
\end{displaymath}
we have
\begin{eqnarray*}
I_{1}&=&\frac{1}{2}\alpha_{2m}(2m-1)!\sum_{i=1}^{k}\int_{B_{\frac{\delta_{0}}{\mu_{i}\varepsilon}}(0)}
\frac{dy}{(1+|y|^{2})^{2m}}\bigg\{-2m\log(\mu_{i}^{2}\varepsilon^{2}+|\mu_{i}\varepsilon y|^{2})\\
&+&H(\xi_{i},\xi_{i})+O(\mu_{i}|\varepsilon y|)+O(\varepsilon^{2})\bigg\}+O(\varepsilon^{2m}).
\end{eqnarray*}
By the fact that $-2m\log(\mu_{i}^{2}\varepsilon^{2}+|\mu_{i}\varepsilon y|^{2})=-2m\log(\mu_{i}^{2}\varepsilon^{2})(1+|y|^{2})$ we have
\begin{eqnarray*}
I_{1}&=&\frac{1}{2}\alpha_{2m}(2m-1)!\sum_{i=1}^{k}\int_{B_{\frac{\delta_{0}}{\mu_{i}\varepsilon}}(0)}
\frac{dy}{(1+|y|^{2})^{2m}}\bigg\{-2m\log(\mu_{i}^{2}\varepsilon^{2})(1+|y|^{2})\\
&+&H(\xi_{i},\xi_{i})+O(\varepsilon |y|)+O(\varepsilon^{2})\bigg\}+O(\varepsilon^{2m}).
\end{eqnarray*}
Now, since
\begin{displaymath}
\int_{B_{\frac{\delta_{0}}{\mu_{i}\varepsilon}}(0)}
\frac{dy}{(1+|y|^{2})^{2m}}O(\varepsilon^{2})=c\ O(\varepsilon^{2}),
\end{displaymath}
because, by improper integral test, the considered integral is bounded and
\begin{displaymath}
\int_{B_{\frac{\delta_{0}}{\mu_{i}\varepsilon}}(0)}
\frac{dy}{(1+|y|^{2})^{2m}}O(\varepsilon |y|)=\int_{B_{\frac{\delta_{0}}{\mu_{i}\varepsilon}}(0)}
\frac{|y| dy}{(1+|y|^{2})^{2m}}O(\varepsilon)=c\ O(\varepsilon)
\end{displaymath}
because, by a comparison test, the considered integral is bounded, thus
\begin{eqnarray*}
I_{1}&=&\frac{1}{2}\alpha_{2m}(2m-1)!\sum_{i=1}^{k}\int_{B_{\frac{\delta_{0}}{\mu_{i}\varepsilon}}(0)}
\frac{dy}{(1+|y|^{2})^{2m}}\bigg\{-2m\log(\mu_{i}^{2}\varepsilon^{2})(1+|y|^{2})\\
&+&H(\xi_{i},\xi_{i})\bigg\}+O(\varepsilon).
\end{eqnarray*}
At this point, we can split the integrals involved in a piece over $\mathbb{R}^{2m}$ and one over $\mathbb{R}^{2m}\setminus B_{\frac{\delta_{0}}{\mu_{i}\varepsilon}}(0)$. By direct computation, we will show that the integrals over $\mathbb{R}^{2m}\setminus B_{\frac{\delta_{0}}{\mu_{i}\varepsilon}}(0)$ are small.
We have
\begin{displaymath}
\int_{\mathbb{R}^{2m}\setminus B_{\frac{\delta_{0}}{\mu_{i}\varepsilon}}(0)}\frac{dy}{(1+|y|^{2})^{2m}}=
2m\omega_{2m}\int_{\frac{\delta_{0}}{\mu_{i}\varepsilon}}^{+\infty}\frac{r^{2m-1}}{r^{4m}}dr=
\frac{(\mu_{i}\varepsilon)^{2m}}{\delta_{0}}=O(\varepsilon^{2m})
\end{displaymath}
and
\begin{displaymath}
\int_{\mathbb{R}^{2m}\setminus B_{\frac{\delta_{0}}{\mu_{i}\varepsilon}}(0)}\frac{\log(1+|y|^{2})dy}{(1+|y|^{2})^{2m}}=
2m\omega_{2m}\int_{\frac{\delta_{0}}{\mu_{i}\varepsilon}}^{+\infty}\frac{r^{2m-1}\log r^{2}}{r^{4m}}dr=
c\int_{\frac{\delta_{0}}{\mu_{i}\varepsilon}}^{+\infty}\frac{\log r}{r^{2m+1}}dr,
\end{displaymath}
integrating by parts and recalling that $\delta_{0}$ is small and we have
\begin{displaymath}
=c\bigg[\log r\int\frac{1}{r^{2m+1}dr}\bigg]_{\frac{\delta_{0}}{\mu_{i}\varepsilon}}^{+\infty}-\int_{\frac{\delta_{0}}{\mu_{i}\varepsilon}}^{+\infty}\frac{1}{r^{2m+1}}\frac{1}{r}dr=O(\varepsilon^{2m}\log\varepsilon).
\end{displaymath}
Now since, clearly, $\varepsilon\leq\varepsilon^{2m}$ and $\varepsilon\leq\varepsilon^{2m}\log\varepsilon$ we have
\begin{eqnarray*}
I_{1}&=&\frac{1}{2}\alpha_{2m}(2m-1)!\sum_{i=1}^{k}\int_{\mathbb{R}^{2m}}
\frac{dy}{(1+|y|^{2})^{2m}}\bigg\{-2m\log(\mu_{i}^{2}\varepsilon^{2})-2m\log(1+|y|^{2})\\
&+&H(\xi_{i},\xi_{i})\bigg\}+O(\varepsilon)=\frac{1}{2}\alpha_{2m}(2m-1)!\bigg[\int_{\mathbb{R}^{2m}}\frac{dy}{(1+|y|^{2})^{2m}}\sum_{i=1}^{k}(H(\xi_{i},\xi_{i})-\\ &&4m\log\mu_{i}\varepsilon)-2m\sum_{i=1}^{k}\int_{\mathbb{R}^{2m}}\frac{\log(1+|y|^{2})dy}{(1+|y|^{2})^{2m}}\bigg]+O(\varepsilon)
\end{eqnarray*}
and, finally,
\begin{equation}\label{E20}
I_{1}=\frac{1}{2}\alpha_{2m}\pi^{m}(m-1)!\sum_{i=1}^{k}(H(\xi_{i},\xi_{i})-4m\log\mu_{i}\varepsilon)-k\alpha_{2m}\pi^{m}(m-1)!+O(\varepsilon),
\end{equation}
where we have used the explicit values (see Appendix)
\begin{displaymath}
c_{0}=\int_{\mathbb{R}^{2m}}\frac{dy}{(1+|y|^{2})^{2m}}=\frac{\pi^{m}(m-1)!}{(2m-1)!}
\end{displaymath}
and
\begin{displaymath}
c_{1}=\int_{\mathbb{R}^{2m}}\frac{\log(1+|y|^{2})dy}{(1+|y|^{2})^{2m}}=\frac{\pi^{m}(m-1)!}{m(2m-1)!}.
\end{displaymath}
Now, we consider the mixed term quadratic part $I_{2}$,
\begin{displaymath}
I_{2}=\frac{1}{2}\rho_{2m}\sum_{i\neq j}V(\xi_{i})\int_{\Omega}e^{u_{i}}U_{j}=\frac{1}{2}\rho_{2m}\sum_{i\neq j}V(\xi_{i})\bigg[\int_{\Omega\setminus B_{\delta_{0}}(\xi_{i})\cup B_{\delta_{0}}(\xi_{j})}e^{u_{i}}U_{j}+
\end{displaymath}
\begin{displaymath}
+\int_{ B_{\delta_{0}}(\xi_{j})}e^{u_{i}}U_{j}+\int_{ B_{\delta_{0}}(\xi_{i})}e^{u_{i}}U_{j}\bigg]=\frac{1}{2}\rho_{2m}\sum_{i\neq j}V(\xi_{i})\bigg[\widetilde{I_{1}}+\widetilde{I_{2}}+\widetilde{I_{3}}\bigg].
\end{displaymath}
Clearly, $\widetilde{I_{1}}=O(\varepsilon^{2m})$. Since $\int_{B_{\delta_{0}}(\xi_{j})}e^{u_{i}}=O(\varepsilon^{2})$ then
\begin{displaymath}
\widetilde{I_{2}}=O(\varepsilon^{2})\int_{B_{\delta_{0}}(\xi_{j})}\log\frac{1}{(\mu_{j}^{2}\varepsilon^{2}+|x-\xi_{j}|^{2})^{2m}}+H(x,\xi_{j})+O(\varepsilon^{2})
\end{displaymath}
\begin{displaymath}
=O(\varepsilon^{2})+O\bigg(\varepsilon^{2}\int_{B_{\delta_{0}}(\xi_{j})}\log\frac{1}{(\mu_{j}^{2}\varepsilon^{2}+|x-\xi_{j}|^{2})^{2m}}\bigg)=O(\varepsilon^{2})
\end{displaymath}
since
\begin{displaymath}
\int_{B_{\delta_{0}}(\xi_{j})}\log(\mu_{j}^{2}\varepsilon^{2}+|x-\xi_{j}|^{2})^{2m}\sim\int_{B_{\delta_{0}}(\xi_{j})}\log(|x-\xi_{j}|)=O(1).
\end{displaymath}
Now,
\begin{displaymath}
\frac{1}{2}\rho^{2m}\sum_{i\neq j}V(\xi_{i})\widetilde{I_{3}}=\frac{1}{2}\rho^{2m}\sum_{i\neq j}V(\xi_{i})\int_{B_{\delta_{0}}(\xi_{i})}\frac{\mu_{i}^{2m}(1+\varepsilon^{2})^{2m}}{(\mu_{i}^{2}\varepsilon^{2}+|x-\xi_{i}|^{2})^{2m}V(\xi_{i})}G(x,\xi_{j})+O(\varepsilon^{2}),
\end{displaymath}
by first order Taylor expansion $G(x,\xi_{j})=G(\xi_{j},\xi_{i})+O(|x-\xi_{j}|)$ and with the same change of variables above we have
\begin{displaymath}
\frac{1}{2}\rho^{2m}\sum_{i\neq j}k(\xi_{i})\widetilde{I_{3}}=\frac{1}{2}\alpha_{2m}(2m-1)!\sum_{i\neq j}\int_{B_{\frac{\delta_{0}}{\mu_{i}\varepsilon}}(0)}\frac{dy}{(1+|y|^{2})^{2m}}G(\xi_{j},\xi_{i})+O(\varepsilon|y|)+O(\varepsilon^{2}).
\end{displaymath}
Since
\begin{displaymath}
\int_{B_{\frac{\delta_{0}}{\mu_{i}\varepsilon}}(0)}\frac{dy}{(1+|y|^{2})^{2m}}=\int_{\mathbb{R}^{2m}}\frac{dy}{(1+|y|^{2})^{2m}}+O(\varepsilon^{2})=c_{0}+O(\varepsilon^{2})
\end{displaymath}
then
\begin{displaymath}
\frac{1}{2}\rho^{2m}\sum_{i\neq j}V(\xi_{i})\widetilde{I_{3}}=\frac{1}{2}\alpha_{2m}\pi^{m}(m-1)!\sum_{i\neq j}G(\xi_{j},\xi_{i})+O\bigg(\varepsilon\int_{B_{\frac{\delta_{0}}{\mu_{i}\varepsilon}}(0)}\frac{|y|}{(1+|y|^{2})^{2m}}\bigg)+O(\varepsilon^{2}).
\end{displaymath}
Since $\displaystyle\int_{B_{\frac{\delta_{0}}{\mu_{i}\varepsilon}}(0)}\frac{|y|}{(1+|y|^{2})^{2m}}=O(1)$ then $O\bigg(\varepsilon\displaystyle\int_{B_{\frac{\delta_{0}}{\mu_{i}\varepsilon}}(0)}\frac{|y|}{(1+|y|^{2})^{2m}}\bigg)=O(\varepsilon)$.\\
Adding all this information we arrive to
\begin{equation}\label{E21}
I_{2}=\frac{1}{2}\alpha_{2m}\pi^{m}(m-1)!\sum_{i\neq j}G(\xi_{j},\xi_{i})+O(\varepsilon).
\end{equation}
Finally, we consider $I_{3}$.
\begin{displaymath}
I_{3}=-\rho^{2m}\int_{\Omega}V(x)e^{U}=-\rho^{2m}\sum_{i=1}^{k}\bigg(\int_{B_{\delta_{0}}(\xi_{i})}V(x)e^{U}+\int_{\Omega\setminus B_{\delta_{0}}(\xi_{i})}V(x)e^{U}\bigg),
\end{displaymath}
since $\int_{\Omega\setminus B_{\delta_{0}}(\xi_{i})}V(x)e^{U}=O(\varepsilon^{2m})$ then
\begin{displaymath}
I_{3}=-\rho^{2m}\sum_{i=1}^{k}\int_{\Omega}V(x)e^{u_{i}(x)+H_{i}(x)}dx+O(\varepsilon^{2m}).
\end{displaymath}
Using the change of variables above we have
\begin{displaymath}
I_{3}=-\rho^{2m}\sum_{i=1}^{k}\int_{B_{\frac{\delta_{0}}{\mu_{i}\varepsilon}}(0)}V(\xi_{i}+\mu_{i}\varepsilon y)\exp\bigg\{u_{i}(\xi_{i}+\mu_{i}\varepsilon y)+H_{i}(\xi_{i}+\mu_{i}\varepsilon y)\bigg\}(\mu_{i}\varepsilon)^{2m}dy+O(\varepsilon^{2m}).
\end{displaymath}
Using \eqref{E18} we have
\begin{eqnarray*}
I_{3}&=&-\rho^{2m}\sum_{i=1}^{k}\int_{B_{\frac{\delta_{0}}{\mu_{i}\varepsilon}}(0)}V(\xi_{i}+\mu_{i}\varepsilon y)\frac{\mu_{i}^{2m}(1+\varepsilon^{2})^{2m}(\mu_{i}\varepsilon)^{2m}}{(\mu_{i}^{2}\varepsilon^{2}+|\mu_{i}\varepsilon y|^{2})^{2m}V(\xi_{i})}\\
&&\exp\bigg\{H(\xi_{i}+\mu_{i}\varepsilon y,\xi_{i})-2m\log\mu_{i}(1+\varepsilon^{2})+\log V(\xi_{i})+O(\varepsilon^{2})\bigg\}dy+O(\varepsilon^{2m}),
\end{eqnarray*}
By a first order Taylor expansion for $H(\xi_{i}+\mu_{i}\varepsilon y,\xi_{i})$ we have
\begin{displaymath}
I_{3}=-\rho^{2m}\sum_{i=1}^{k}\int_{B_{\frac{\delta_{0}}{\mu_{i}\varepsilon}}(0)}\frac{V(\xi_{i}+\mu_{i}\varepsilon y)}{(\mu_{i}\varepsilon)^{2m}(1+|y|^{2})^{2m}}\cdot
\end{displaymath}
\begin{displaymath}
\exp\bigg\{H(\xi_{i},\xi_{i})+O(|\mu_{i}\varepsilon y|)+O(\varepsilon^{2})\bigg\}dy+O(\varepsilon^{2m})=
\end{displaymath}
\begin{displaymath}
-\alpha_{2m}(2m-1)!\frac{\varepsilon^{2m}}{(1+\varepsilon^{2})^{2m}}\sum_{i=1}^{k}\int_{B_{\frac{\delta_{0}}{\mu_{i}\varepsilon}}(0)} \frac{V(\xi_{i}+\mu_{i}\varepsilon y)}{(\mu_{i}\varepsilon)^{2m}(1+|y|^{2})^{2m}}\cdot
\end{displaymath}
\begin{displaymath}
\exp\bigg\{H(\xi_{i},\xi_{i})+O(\varepsilon |y|)+O(\varepsilon^{2})\bigg\}dy+O(\varepsilon^{2m}),
\end{displaymath}
since $(1+\varepsilon^{2})^{2m}\sim 1+O(\varepsilon^{2m})$ then
we have
\begin{displaymath}
I_{3}=-\alpha_{2m}(2m-1)!\sum_{i=1}^{k}\int_{B_{\frac{\delta_{0}}{\mu_{i}\varepsilon}}(0)}\frac{V(\xi_{i}+\mu_{i}\varepsilon y)}{\mu_{i}^{2m}(1+|y|^{2})^{2m}}e^{H(\xi_{i},\xi_{i})}(1+O(\varepsilon|y|))(1+O(\varepsilon^{2}))dy+O(\varepsilon^{2m}).
\end{displaymath}
By first order Taylor expansion for $V(\xi_{i}+\mu_{i}\varepsilon y)$, the fact that $(1+O(\varepsilon|y|))(1+O(\varepsilon^{2}))=1+O(\varepsilon^{2})+O(\varepsilon|y|)+O(\varepsilon^{2}+\varepsilon|y|)\sim 1+O(\varepsilon)$, the election of $\mu_{i}$'s \eqref{E15} and the integral value $c_{0}$ we have, finally,
\begin{equation}\label{E22}
I_{3}=-\alpha_{2m}\pi^{m}(m-1)!k+O(\varepsilon).
\end{equation}
Thus, summing up \eqref{E20}, \eqref{E21}, \eqref{E22}, we can conclude the following expansion of $J_{\rho}(U)$
\begin{eqnarray*}
J_{\rho}(U)&=&\underbrace{\frac{1}{2}\alpha_{2m}\pi^{m}(m-1)!\bigg[\sum_{i=1}^{k}H(\xi_{i},\xi_{i})-4m\log\mu_{i}\varepsilon+\sum_{i\neq j}G(\xi_{j},\xi_{i})\bigg]}_{I_{4}}\\
&&-2k\alpha_{2m}\pi^{m}(m-1)!.
\end{eqnarray*}
We focus on $I_{4}$.
\begin{displaymath}
I_{4}=\frac{1}{2}\alpha_{2m}\pi^{m}(m-1)!\bigg[\sum_{i=1}^{k}H(\xi_{i},\xi_{i})-4m\log\mu_{i}-4m\log\varepsilon+\sum_{i\neq j}G(\xi_{j},\xi_{i})\bigg],
\end{displaymath}
by the election of $\mu_{i}$'s \eqref{E15} we have
\begin{eqnarray*}
I_{4}&=&\frac{1}{2}\alpha_{2m}\pi^{m}(m-1)!\bigg[\sum_{i=1}^{k}H(\xi_{i},\xi_{i})-2\log V(\xi_{i})-2H(\xi_{i},\xi_{i})-2\sum_{i\neq j}G(\xi_{j},\xi_{i})\bigg]\\
&+&\frac{1}{2}\alpha_{2m}\pi^{m}(m-1)!\sum_{i\neq j}G(\xi_{j},\xi_{i})-2m\alpha_{2m}\pi^{m}(m-1)!\sum_{i=1}^{k}|\log\varepsilon|\\
&=&\frac{1}{2}\alpha_{2m}\pi^{m}(m-1)!\bigg[\sum_{i=1}^{k}-H(\xi_{i},\xi_{i})-2\log V(\xi_{i})\bigg]\\
&-&\frac{1}{2}\alpha_{2m}\pi^{m}(m-1)!\sum_{i\neq j}G(\xi_{j},\xi_{i})-2mk\alpha_{2m}\pi^{m}(m-1)!|\log\varepsilon|.
\end{eqnarray*}
Introducing the function $\varphi_{k}(\xi)$ defined as in \eqref{E17} and the new constant $b_{m}$ defined as in the statement we conclude the proof.
\begin{flushright}
$\Box$
\end{flushright}
In the following, we will stay in the expanded variable $y\in\Omega_{\varepsilon}$.\\ We will look for solutions of problem \eqref{E1} in the form of a small perturbation of $V$, the first approximation of \eqref{E13}. Since we have, by Lemma \ref{b}, a small $\emph{error term}$ $R$, the equation for the perturbation is a linear one with right hand side given by this error term $R$ perturbed by a lower order nonlinear term. The mapping properties of this linear operator are fundamental in solving for such a perturbation. Not only this, the nonlinearity must remain small if, say, an iterative scheme is produced. An obvious way to write this perturbation is in additive way, say $w=W+\phi$ where the $\emph{remainder term}$ $\phi$ will represent a lower order correction, that is $\phi$ goes to zero as $\rho$ goes to zero. The $\emph{nonlinearity term}$ $N(\phi)$ produced when substituting this new ansatz $w$ in \eqref{E1} is a polynomial in $\phi$ carrying at least quadratic terms.
We aim at finding solutions for $\phi$ small provided that the points $\xi_{i}$ are suitably chosen.\\ In term of a small $\phi$, we can rewrite problem \eqref{E1} as a nonlinear perturbation of its linearization $\mathcal{L}_{\varepsilon}(\phi):=(-\Delta)^{m}\phi -T(y)\phi$, namely,
\begin{equation}\label{E23}
\left\{ \begin{array}{ll}
         \mathcal{L}_{\varepsilon}(\phi) =-R+N(\phi)& \mbox{in}\ \Omega_{\varepsilon}\\
          B_{j}\phi=0, |j|\leq m-1& \mbox{on}\ \partial\Omega_{\varepsilon},\
        \end{array} \right.
\vspace{0,2cm}\end{equation}
where
\begin{equation}\label{E24}
T(y):=V(\varepsilon y)e^{W(y)},
\end{equation}
\begin{equation}\label{E25}
R(y):=(-\Delta)^{m}W-T(y),
\end{equation}
\begin{equation}\label{E26}
N(\phi):=T(y)(e^{\phi}-\phi-1).
\end{equation}
Let us observe that near the concentration points $\xi_{i}$ the linearized operator $\mathcal{L}_{\varepsilon}$ is a small nontrivial perturbation of the polyharmonic operator while it is essentially this operator in most of the domain.\\
Now we intend to solve \eqref{E23}. To do so, we need to analyze the possibility to invert the operator $\mathcal{L}_{\varepsilon}(\phi)$ in order to express the equation as a fixed point problem. It is not expected this operator to be, in general, globally invertible because, when regarded in the entire $\mathbb{R}^{2m}$, this operator does have kernel: functions $Y_{i0}$ and $Y_{ij}$, defined in \eqref{E28} and \eqref{E29} below, with $j=1,...,2m$, $i=1,...,k$, annihilate it.
In suitable spaces, for instance $L^{\infty}$, these functions span the entire kernel. In a suitable orthogonal to this kernel, the bilinear form associated to this operator turns out to be uniformly positive definite. Then we are intended to solve a suitably projected version of our problem for which a linear theory is in order and, after which, the resolution comes from a direct application of contraction mapping principle. The next step will be adjust the points $\xi$ in order to have solutions to the full problem. The latter problem will turn out to be equivalent to a variational problem in $\xi$. Our main results will be a consequence of solving this finite dimensional problem. Moreover, we have the validity of the following estimates
\begin{lm}\label{d}
\end{lm}
\emph{For $y\in\Omega_{\varepsilon}$ we have that
\begin{enumerate}
  \item $\|T\|_{\ast} = O(1)$,
  \item $\|N(\phi_{1}) - N(\phi_{2})\|_{\ast} \leq C\max_{i=1,2}\|\phi_{i}\|_{\infty}\|\phi_{1}-\phi_{2}\|_{\infty}$ ($N$ is a contraction),
  \item $\|N(\phi)\|_{\ast} = O(\|\phi\|_{\infty}^{2})$ \mbox{as $\|\phi\|_{\infty} \rightarrow 0$} ($N$ is almost quadratic).
\end{enumerate}}
$\underline{Proof}$\\
From Lemma \ref{a} we have
\begin{displaymath}
V(\varepsilon y)e^{w(y)}=\alpha_{2m}(2m-1)!\frac{\mu_{i}^{2m}}{(\mu_{i}^{2}+|y-\xi_{i}'|^{2})^{2m}}(1+O(\varepsilon |y-\xi_{i}'|))
\end{displaymath}
that is
\begin{displaymath}
T(y)=\alpha_{2m}(2m-1)!T_{0}(y)(1+O(\varepsilon |y-\xi_{i}'|))
\end{displaymath}
where $T_{0}(y):=\sum_{i=1}^{k}\frac{\mu_{i}^{2m}}{(\mu_{i}^{2}+|y-\xi_{i}'|^{2})^{2m}}$. This fact gives the validity of the first claim.\\
For the general problem $(-\Delta)^{m}u=V(x)f(u)$ in a domain $\Omega$ it is possible to rewrite \eqref{E25} and \eqref{E26} in the following way, respectively:
\begin{eqnarray*}
R(y)&=&(-\Delta)^{m}W+V(x)f(W)\\
N(\phi)&=&V(x)[f(W+\phi)-f(W)-f'(W)\phi]
\end{eqnarray*}
where $W$ is the first approximation of the predicted solution and $f'$ is a functional derivative in the Frech\'{e}t sense. Hence, from this, we gain that $\mathcal{L}(\phi)=(-\Delta)^{m}\phi+V(x)f'(W)\phi$ is the related linearized operator.\\
We calculate
\begin{eqnarray*}
|N(\phi_{1})-N(\phi_{2})|&=&V(x)|f(W+\phi_{1})-f(W+\phi_{2})-f'(W)(\phi_{1}-\phi_{2})|\\
&=&V(x)\bigg|\int_{0}^{1}f'(W+t\phi_{1}+(1-t)\phi_{2})dt- \int_{0}^{1}f'(W)(\phi_{1}-\phi_{2})dt\bigg|\\
&=&V(x)\bigg|\int_{0}^{1}(\phi_{1}-\phi_{2})dt + \int_{0}^{1}f''(W+st\phi_{1}+s(1-t)\phi_{2})\cdot\\
&&(t\phi_{1}+(1+t)\phi_{2}))ds\bigg|,
\end{eqnarray*}
and by Lagrange's theorem we have
\begin{displaymath}
\leq \sup |f''||\phi_{1}-\phi_{2}|\int_{0}^{1}|t\phi_{1}-(1-t)\phi_{2}|dt \leq \sup |f''||\phi_{1}-\phi_{2}| \max (|\phi_{1}|;|\phi_{2}|).
\end{displaymath}
This fact gives the validity of the second claim.\\
Now, using this second claim with $\phi_{1}=\phi$ and $\phi_{2}=0$ we, easily, get the third.
\begin{flushright}
$\Box$
\end{flushright}
\begin{oss}\label{R1}
\end{oss}
\begin{enumerate}
  \item For the previous analysis, using $2.$, $3.$ of Lemma \ref{d}, we can prove that also $N-R$ is a contraction. In fact,
      $$
      \|(N(\phi_{1}) - R) - (N(\phi_{2}) - R) \|_{\ast}=\|N(\phi_{1}) - N(\phi_{2})\|_{\ast}\leq C\max_{i=1,2}\|\phi_{i}\|_{\infty}\|\phi_{1}-\phi_{2}\|_{\infty}
      $$
      and, recalling that $R$ is small by Lemma \ref{b},
       \begin{displaymath}
      \|(N(\phi) - R)\|_{\ast}\leq \|N(\phi)\|_{\ast}+ \|R\|_{\ast}\leq C\|\phi\|_{\infty}^{2}+O(\varepsilon).
      \end{displaymath}
  \item We note that far from the points $\xi_{i}'$ (i.e. on most of the domain $\Omega$) $T(y) = O(\varepsilon^{2n}$), hence $\mathcal{L}_{\varepsilon}(\phi)$ is a \emph{small perturbation} of $(-\Delta)^{m}$ away for the concentration points while near the points $\xi_{i}'$, $T(y) = O(\sum_{i=1}^{k}e^{w}_{i}(y))$, hence $\mathcal{L}_{\varepsilon}(\phi)$ is a \emph{nontrivial perturbation} of $(-\Delta)^{m}$ thus $\mathcal{L}_{\varepsilon}(\phi)$ is approximately a superposition of the linear operators $\mathcal{L}_{i}(\phi)$ = $(-\Delta)^{m}\phi - \sum_{i=1}^{k}e^{w_{i}}\phi$.
\end{enumerate}
\section{Analysis of the linearized operator}\label{an}
A main step in solving \eqref{E23} for small $\phi$ under a suitable choice of the point $\xi_{i}$ is that of a solvability theory for the $2m$-order linear operator $\mathcal{L}_{\varepsilon}$ under suitable orthogonality conditions. We shall devote this section to prove the bounded invertibility of $\mathcal{L}_{\varepsilon}$. We consider in $\Omega_{\varepsilon}$
\begin{displaymath}
\mathcal{L}_{\varepsilon}(\phi):=(-\Delta)^{m}\phi - T(y)\phi
\end{displaymath}
for function $\phi$ defined on $\Omega_{\varepsilon}$ and where $T(y)$ was introduced in \eqref{E24}.\\
Unlike the $2m$- Laplacian $(-\Delta)^{m}$, the operator $\mathcal{L}_{\varepsilon}$ as an approximate kernel which in principle prevents any form of bounded invertibility. In fact, centering the system of coordinates at, say, $\xi'_{i}$ by setting $z=y-\xi'_{i}$, one can easily see that formally $\mathcal{L}_{\varepsilon}$ as $\varepsilon\rightarrow 0$ can be approximately regarded as a superposition of the linear operators in $\mathbb{R}^{2m}$
\begin{equation}\label{E27}
\mathcal{L}_{i}(\phi):=(-\Delta)^{m}\phi -\sum_{i=1}^{k}e^{w_{i}}\phi=(-\Delta)^{m}\phi -\alpha_{2m}(2m-1)!\sum_{i=1}^{k}\frac{\mu_{i}^{2m}}{(\mu_{i}^{2}+|z|^{2})^{2m}}\phi,
\end{equation}
namely, Liouville equation $(-\Delta)^{m}w=e^{w}$ linearized around the standard bubble $$w_{i}(z)=\alpha_{2m}(2m-1)!\sum_{i=1}^{k}\log\frac{\mu_{i}^{2m}}{(\mu_{i}^{2}+|z|^{2})^{2m}}.$$\\
Thus the key point to develop a satisfactory solvability theory for the operator $\mathcal{L}_{\varepsilon}$ is the $\emph{non-degeneracy}$ of the solutions $w_{i}$ up to the natural invariances of the Liouville equation under translations and, due to the presence of the critical exponent, also under dilations, that is if $z$ is a solution, necessarily also the function $x\mapsto\mu^{\alpha}z(\frac{x-\xi}{\mu})$ is a solution, for any $\mu>0$ and for some suitable constant $\alpha$.\\ Non-degeneracy is an important ingredient in the construction of solutions to problems involving small parameters and concentration phenomena in which, after a suitable blowing-up around a concentration point, one sees a limiting equation. This property, in general, is used to build solutions with multiple concentration points.
We are interested in the classification of bounded solutions of $\mathcal{L}_{\varepsilon}(\phi)=0$ in $\mathbb{R}^{2m}$. Some bounded solutions are easy to find. For example, we can define the functions
\begin{equation}\label{E28}
Y_{i0}(z)= -\frac{r}{2m}\partial_{r}(v_{i})+2m = \frac{-\mu_{i}^{2}+|z|^{2}}{\mu_{i}^{2}+|z|^{2}},
\end{equation}
where $r=|z|$ and, clearly, $\mathcal{L}_{i}(Y_{i0})=0$ and this reflect the fact that the equation is invariant under the group of dilations $\tau\rightarrow -u (\tau \cdot) + 2m \log \tau$ and for all $j=1,...,2m$ and $i=1,...,k$
\begin{equation}\label{E29}
Y_{ij}(z)= -\partial_{z_{i}} v_{i} = 4m\frac{z_{i}}{\mu_{i}^{2}+|z|^{2}},
\end{equation}
clearly $\mathcal{L}_{i}(Y_{ij})=0$ and this reflect the fact that the equation is invariant under the group of translations $a\rightarrow u (\cdot +a)$.\\
In the next Lemma, following very close the geometric point of view in (\cite{barak2}, Lemma 3.1), we prove that any solutions of the linearized operator $\mathcal{L}_{i}$ calculated along the standard bubbles are nondegenerate, that is it has not trivial kernel in certain space, for instance $L^{\infty}$, and thus, by invariances of the Liouville equation under a large groups of symmetries, $\mathcal{L}_{i}$ has bounded kernel. This is equivalent to the fact that the actually $L^{\infty}$-kernel of the linearized operator $\ker\mathcal{L}_{i}\cap L^{\infty}(\mathbb{R}^{2m})$ is spanned by the bounded solutions \eqref{E28} and \eqref{E29} which naturally belong to this space, that is the only bounded solution of $\mathcal{L}_{i}(\phi)=0$ in all $\mathbb{R}^{2m}$ are linear combinations of those functions. In the proof, we use some well-known fact in spectral theory and the classical idea to prove nondegeneracy comparing precisely the dimension of the \emph{bubble solutions manifold} with respect to the dimension of the kernel of the linearized operator calculated along the standard bubble. As a byproduct, we have that $\ker \mathcal{L}_{i}$ is equal to the eigenspace of $\mathbf{P}_{g_{\mathbb{S}^{2m}}}$ associated to the eigenvalue $t_{m}$ and that $\mbox{dim} (\ker \mathcal{L}_{i}) = 2m+1$.
\begin{lm}\label{e}($L^{\infty}$ non-degeneracy)
\end{lm}
\emph{Any bounded solutions of $\mathcal{L}_{i}(\phi)=0$ in all $\mathbb{R}^{2m}$ are linear combinations of $\eqref{E28}$ and $\eqref{E29}$.
}\\$\underline{Proof}$\\
We consider on $\mathbb{R}^{2m}$ the Euclidean metric $g_{E}=|dz|^{2}$ and the spherical metric on $\mathbb{S}^{2m}$
\begin{displaymath}
g_{\mathbb{S}^{2m}} = \frac{4}{(1+|z|^{2})^{2}}g_{E}
\end{displaymath}
induced by the inverse of the stereographic projection from the sphere $\mathbb{S}^{2m}$ onto the whole $\mathbb{R}^{2m}$ with respect to the North pole, namely
\begin{displaymath}
\pi^{-1}: z=(z_{j}) \in \mathbb{R}^{2m} \rightarrow
\mathbb{S}^{2m}\setminus \{\mbox{North\ pole}\}\subset \mathbb{R}^{2m+1}
\end{displaymath}
given by the formulas
\begin{displaymath}
\pi^{-1}(z):=
\begin{cases}
\frac{2z_{j}}{1+|z|^{2}} & 1\leq j\leq 2m\\
\frac{1-|z|^{2}}{1+|z|^{2}} & j= 2m+1.
\end{cases}
\end{displaymath}
We remark that $\pi^{-1}$ is a conformal diffeomorphism, more precisely the pullback of $g_{\mathbb{S}^{2m}}$ to $\mathbb{R}^{2m}$ satisfies
\begin{displaymath}
(\pi^{-1})^{\ast}g_{\mathbb{S}^{2m}} = \frac{4}{(1+|z|^{2})^{2}}g_{E}.
\end{displaymath}
According to \cite{beck,bran2}, we have $\mathbf{P}_{g_{\mathbb{S}^{2m}}} = \prod_{k=0}^{m-1}(-\Delta_{\mathbb{S}^{2m}}+k(2m-k-1))$ where $\mathbf{P}_{g_{\mathbb{S}^{2m}}}$ is the \emph{GJMS}-operator on $\mathbb{S}^{2m}$ and $-\Delta_{\mathbb{S}^{2m}}$ is the Laplace-Beltrami operator on $\mathbb{S}^{2m}$ and when the manifold is the Euclidean space, the \emph{GJMS}-operator is simply given by $\mathbf{P}_{\mathbb{R}^{2m}} = (-\Delta)^{m}$. Since under the following conformal change of metric $\widetilde{g} = e^{2\phi}g$ the \emph{GJMS}-operator transforms according to $\mathbf{P}_{\widetilde{g}} = e^{-2m\phi}\mathbf{P}_{g}\phi$ we obtain that $\bigg(\frac{4}{(1+|z|^{2})^{2}}\bigg)^{m}\mathbf{P}_{g_{\mathbb{S}^{2m}}} = \mathbf{P}_{g_{E}}$.\\
In particular, if $\phi \equiv \phi\circ \pi: (\mathbb{R}^{2m}, g_{E})\rightarrow \mathbb{R}$ (with slight abuse of notation we identify $\phi$ with $\phi\circ \pi$) is a bounded solution of $\mathcal{L}_{i}(\phi)=0$ then $\phi \equiv \phi\circ \pi: (\mathbb{R}^{2m}, g_{\mathbb{S}^{2m}})\rightarrow \mathbb{R}$ is a bounded solution of
\begin{equation}\label{E30}
\mathbf{P}_{g_{\mathbb{S}^{2m}}}\phi = t_{m}\phi\ \mbox{in}\ \mathbb{S}^{2m}\setminus \{\mbox{North\ pole}\}
\end{equation}
away from the North pole and where $t_{m}:= \frac{\alpha_{2m}(2m-1)!}{2^{2m}}=2m(m-1)!(2m-1)!$. For instance, $t_{2}=2, t_{4}=24$.\\
Since $\phi$ is assumed to be bounded then the isolated singularity at the North pole is removable and hence \eqref{E30} holds on all $\mathbb{S}^{2m}$.\\ We now perform the eigenfunction decomposition of $\phi$ in terms of the eigendata of the Laplacian on $\mathbb{S}^{2m}$. We decompose
\begin{displaymath}
\phi=\sum_{k\geq0}\phi_{k}
\end{displaymath}
where $\phi_{k}$ belongs to the $k$-th eigenspace of $\Delta_{S^{2m}}$, namely $\phi_{k}$ satisfies the following eigenvalue's problem $-\Delta_{g_{\mathbb{S}^{2m}}}\phi_{k} = \lambda_{k}\phi_{k}$. It is well known that the spectrum of $-\Delta_{g_{\mathbb{S}^{2m}}}$ is discrete and represented by
\begin{displaymath}
\sigma(-\Delta_{g_{\mathbb{S}^{2m}}}) = \{\lambda_{k} /k\geq 0\} = \{k(k+2m-1) /k\geq 0\}
\end{displaymath}
with multiplicity $m_{k}=\frac{(2m+k-2)!(2m+2k-1)!}{k!(2m-1)!}$, see (\cite{dupa}, Appendix C).
By (\ref{E30}) we have
\begin{displaymath}
\prod_{k=0}^{m-1}\bigg[(-\lambda_{k}+k(2m-k-1))-t_{m}\bigg]\phi_{k}=0,
\end{displaymath}
that is $\prod_{k=0}^{m-1}(-2k^{2}-t_{m})\phi_{k}=0$. Hence $\phi_{k}=0$ for all $k\geq0$ except, eventually, those for which $\prod_{k=0}^{m-1}(-2k^{2}-t_{m})=0$. This implies that $\phi: \mathbb{S}^{2m}\rightarrow \mathbb{R}$ is a combination of the eigenfunctions associated to $k=1$ that are given by $\varphi_{j}(y) = y_{j}$, $j=1,...,2m+1$, where $y=(y_{j})\in \mathbb{S}^{2m}$. Being the sphere parameterized by $\pi^{-1}$ we may write $y=\pi^{-1}(z)$. Then, the functions $2m\varphi_{j}$ precisely correspond to the functions $Y_{ij}$ while the function $2m\varphi_{2m+1}$ corresponds to the function $Y_{i0}$.
\begin{flushright}
$\Box$
\end{flushright}
\section{Projected linear theory for $\mathcal{L}_{\varepsilon}$ onto kernel}\label{pro}
We define for $j=0,...,2m$ and $i=1,...,k$,
\begin{displaymath}
Z_{ij}(y):=Y_{ij}(z).
\end{displaymath}
Additionally, let us consider $R_{0}$ a large but fixed number and $\chi:\mathbb{R}\rightarrow \mathbb{R}$, $\chi_{i}(y)=\chi(r_{i})$, where $r_{i}:=|z|$, a smooth and radial cut-off function with $\chi:=1$ if $r_{i}\leq R_{0}$, namely in $B_{R_{0}}(0)$ and $\chi:=0$ if $r_{i}\geq R_{0}+1$, namely in the complementary $B_{R_{0}+1}(0)^{c}$.\\
Given $h\in L^{\infty}(\Omega_{\varepsilon})$, we consider the linear problem of finding a function $\phi:\Omega_{\varepsilon}\rightarrow \mathbb{R}$, $\phi=\phi(y)$, such that for certain scalars $c_{ij}$ one has
\begin{equation}\label{E31}
\left\{ \begin{array}{ll}
         \mathcal{L}_{\varepsilon}(\phi)=h+\sum_{j=1}^{2m}\sum_{i=1}^{k}c_{ij}\chi_{i}Z_{ij}& \mbox{in}\ \Omega_{\varepsilon},\\
          B_{j}\phi=0, |j|\leq m-1& \mbox{on}\ \partial\Omega_{\varepsilon},\\
          \int_{\Omega_{\varepsilon}}\chi_{i}Z_{ij}\phi=0& \mbox{for all}\ j=1,...,2m,\ i=1,...,k.
        \end{array} \right.
\vspace{0,2cm}\end{equation}
The orthogonality condition in \eqref{E31} are only taken with respect to the elements of the approximate kernel due to translations.\\
The main goal of this section is the bounded solvability of \eqref{E31}. Before this, we will establish a priori estimates for this problem. To this end we shall conveniently introduce an adapted norm in $\Omega_{\varepsilon}$. Given $\phi:\Omega_{\varepsilon}\rightarrow \mathbb{R}$ and $\alpha\in \mathbb{N}^{k}$ a multi-index of order $k$ and length $|\alpha|$ we define
\begin{equation}
\|\phi\|_{\ast\ast}:=\sum_{i=1}^{k}\|\phi\|_{C^{2m,\alpha}(r_{i}<2)}+\sum_{i=1}^{k}\sum_{|\alpha|\leq 2m-1}\|r_{i}^{|\alpha|}D^{\alpha}\phi\|_{L^{\infty}(r_{i}\geq 2)}.
\end{equation}
We remark that the \emph{interior portion} of the norm defined above controls the function $\phi$ in a neighborhood of the origin while the \emph{exterior portion} the decay at infinity of $\phi$.
With this definition at hand, we prove the main result of this section:
\begin{pro}\label{P1}
\end{pro}
\emph{There exist positive constants $\varepsilon_{0}$ and $C$ such that for any $h\in L^{\infty}(\Omega_{\varepsilon})$, with $\|h\|_{\ast}<\infty$, and any $\xi\in\mathcal{M}$, there is a unique solution $\phi=Q(h)$ to problem \eqref{E31} for all $\varepsilon$ sufficiently small, say $\varepsilon\in(0,\varepsilon_{0})$, which defines a linear operator of $h$. Besides, we have the a priori estimate
\begin{equation}\label{Eb}
\|\phi\|_{\ast\ast}\leq C|\log\varepsilon|\|h\|_{\ast}.
\end{equation}}\\
Before proceeding with the proof we remember a useful lemma on the removable singularities for polyharmonic function. Denote by $B$ and $B_{0}$ the $n$-unit ball and the $n$-punctured unit ball $B\setminus\{0\}$, respectively. We present the result in (\cite{futa}, Theorem 3.2) specified to the $n=2m$ case.
\begin{thm}\label{T2}
\end{thm}
\emph{Suppose $n=2m$ and $u$ a $m$-harmonic function defined on $B_{0}$. Then, the following are equivalent:
\begin{enumerate}
  \item $u$ can be extended to a m-harmonic function on $B$;
  \item $\displaystyle\lim_{x\rightarrow0}u(x)$ exists and is finite;
  \item $u$ is bounded near the origin.\\
\end{enumerate}}

The proof of Proposition \ref{P1} will be split into a series of lemmas which we state and prove next. The first step is to obtain a priori estimates for the problem
\begin{equation}\label{E32}
\left\{ \begin{array}{ll}
         \mathcal{L}_{\varepsilon}(\phi)=h& \mbox{in}\ \Omega_{\varepsilon},\\
          B_{j}\phi=0, |j|\leq m-1& \mbox{on}\ \partial\Omega_{\varepsilon},\\
          \int_{\Omega_{\varepsilon}}\chi_{i}Z_{ij}\phi=0& \mbox{for all}\ j=0,...,2m,\ i=1,...,k
        \end{array} \right.
\vspace{0,2cm}\end{equation}
which involves more orthogonality conditions than those in \eqref{E31}. Notice that in the case $m=1$, independently from the nonlinearity, a key step in order to prove such result is the fact that the operator $\mathcal{L}_{\varepsilon}$ satisfies maximum principle in $\Omega_{\varepsilon}$ outside large balls, see (\cite{delpi}, Lemma 3.1) and \cite{esp} for the exponential-type nonlinearity or the surveys \cite{pis} and references therein for other nonlinearity issue such as the Brezis-Nirenberg Problem or the Coron's Problem. For our $\mathcal{L}_{\varepsilon}$ this is not more true and we need a different approach.
We have the following estimate
\begin{lm}\label{aa}
\end{lm}
\emph{There exist positive constants $\varepsilon_{0}$ and $C$ such that for any solution $\phi$ of problem \eqref{E32} with $h\in L^{\infty}(\Omega_{\varepsilon})$, $\|h\|_{\ast}<\infty$ and with $\xi\in\mathcal{M}$, then for all $\varepsilon\in(0,\varepsilon_{0})$
\begin{equation}\label{E33}
\|\phi\|_{\ast\ast}\leq C\|h\|_{\ast}.
\end{equation}
}\\$\underline{Proof}$\\
We carry out the proof by means of a contradiction argument. If the above fact were false, then, suppose there exist a sequence $\varepsilon_{n}\rightarrow0$, a sequence of points $\xi^{n}=(\xi^{n}_{i})\in\mathcal{M}$ with $i=1,...,k$, a sequence of functions $h_{n}$ with $\|h_{n}\|_{\ast}\rightarrow0$ as $n\rightarrow\infty$ and a sequence of associated solutions $\phi_{n}$ with $\|\phi_{n}\|_{\ast\ast}=1$ of the following problem
\begin{equation}
\left\{ \begin{array}{ll}
         \mathcal{L}_{\varepsilon_{n}}(\phi_{n})=h_{n}& \mbox{in}\ \Omega_{\varepsilon_{n}},\\
          B_{j}\phi_{n}=0, |j|\leq m-1& \mbox{on}\ \partial\Omega_{\varepsilon_{n}},\\
          \int_{\Omega_{\varepsilon_{n}}}\chi_{i}Z_{ij}\phi_{n}=0& \mbox{for all}\ j=0,...,2m,\ i=1,...,k,
        \end{array} \right.
\vspace{0,2cm}\end{equation}
where $\mathcal{L}_{\varepsilon_{n}}(\phi_{n}):=(-\Delta)^{m}\phi_{n}-V(\varepsilon_{n}y)e^{W(y)}\phi_{n}$.\\
We observe that $\|h_{n}\|_{\ast}\rightarrow0$ as $n\rightarrow\infty$ because if \eqref{E33} doesn't hold then $\|\phi_{n}\|_{\ast\ast}> C\|h_{n}\|_{\ast}$, normalizing this last inequality dividing by $\|\phi_{n}\|_{\ast\ast}$ we have done.
We will show that, under those hypothesis, we necessarily have $\phi_{n}\rightarrow0$ in $\Omega_{\varepsilon}$ and, then, a contradiction arises. \\Let us set $\tilde{\phi}_{n}:=\phi_{n}(\varepsilon_{n}^{-1}x)$, $x\in\Omega$.
First, we prove an intermediate claim: with the above hypothesis we have that $\tilde{\phi}_{n}$ goes to zero as $n\rightarrow\infty$ in $C^{m+1,\alpha}$-sense uniformly over compact subsets of $\Omega\setminus \{\xi^{\ast}\}$. In particular, for any $\delta'>0$ sufficiently small we have $$\|\tilde{\phi}_{n}\|_{L^{\infty}(\Omega\setminus\bigcup_{i=1}^{k}B_{\delta'}(\xi_{n})}\rightarrow0\ \mbox{as}\ n\rightarrow\infty.$$
Now, $\tilde{\phi}_{n}$ solves the equation in $\Omega_{\varepsilon_{n}}$:
\begin{displaymath}
(-\Delta)^{m}\tilde{\phi}_{n}(x)=\varepsilon_{n}^{-2m}(-\Delta)^{m}\phi_{n}(\varepsilon_{n}^{-1}x)=\varepsilon_{n}^{-2m}[-T_{n}(\varepsilon_{n}^{-1}x)\phi_{n}(\varepsilon_{n}^{-1}x)+h_{n}(\varepsilon_{n}^{-1}x)]
\end{displaymath}
\begin{displaymath}
=\varepsilon_{n}^{-2m}[V(x)e^{W(\varepsilon_{n}^{-1}x)}\tilde{\phi}_{n}(x)-h_{n}(\varepsilon_{n}^{-1}x)]=\varepsilon_{n}^{-2m}[V(x)e^{U(x)}\tilde{\phi}_{n}(x)-h_{n}(\varepsilon_{n}^{-1}x)].
\end{displaymath}
For any $\delta'>0$ sufficiently small, we recall that in $|y-(\xi^{n})'|>\varepsilon_{n}^{-1}\delta'$ where $(\xi^{n})':=\varepsilon_{n}^{-1}\xi^{n}$ we have $\varepsilon_{n}^{2m}T_{n}(y)\phi_{n}(y)=O(\varepsilon_{n}^{2m})$ because in this region $T_{n}(y)=O(\varepsilon_{n}^{2m})$, that is small, and $\phi_{n}$ is uniformly bounded by contradiction hypothesis. Setting $x=\varepsilon_{n}y$ and recalling that $W(y)=U(x)+2m\log\rho\varepsilon_{n}$, then also $\varepsilon_{n}^{2m}V(x)e^{U(x)}\tilde{\phi_{n}}(x)=O(\varepsilon_{n}^{2m})$. Thus, $(-\Delta)^{m}\tilde{\phi}_{n}(x)=O(\varepsilon_{n}^{2m})+\varepsilon_{n}^{-2m}h_{n}(\varepsilon_{n}^{-1}x)$. \\At this point, we claim that $O(\varepsilon_{n}^{2m})+\varepsilon_{n}^{-2m}h_{n}(\varepsilon_{n}^{-1}x)=o(1)$ as $n\rightarrow\infty$ uniformly in $\widetilde{\Omega}:=\Omega\setminus\bigcup_{i=1}^{k}B_{\delta'}(\xi_{n})$. In fact, if $|x-\xi_{n}|>\delta'$ then the dominant part in the definition of $\|\cdot\|_{\ast}$ is, as noted before, the second, thus, $\varepsilon_{n}^{-2m}|h_{n}(y)|\leq c\|h_{n}(x)\|_{\ast}$ and by contradiction hypothesis we easily get the desired result. We have directly checked that for any $\delta'>0$ small $\tilde{\phi}_{n}$ solves the problem
$$
\left\{ \begin{array}{ll}
         (-\Delta)^{m}\tilde{\phi}_{n}=o(1)& \mbox{uniformly in}\ \widetilde{\Omega},\\
          B_{j}\tilde{\phi}_{n}=0, |j|\leq m-1& \mbox{on}\ \partial\Omega,\
        \end{array} \right.
\vspace{0,2cm}$$
together with the sharp estimates:
\begin{equation}\label{E34}
\|\tilde{\phi}_{n}\|_{L^{\infty}(\widetilde{\Omega)}}\leq1\ \mbox{and}\  \|B_{\tilde{j}}\tilde{\phi}_{n}\|_{L^{\infty}(\widetilde{\Omega)}}\leq C_{\delta'},\ \tilde{j}=1,...,m-1.
\end{equation}
Differently from (\cite{clapp}, Lemma 3.1) we are also able to compute the constant $C_{\delta'}$.
\\Fix an index $0\leq q\leq 2m-1$, since by contradiction hypothesis $\|\phi_{n}\|_{\ast\ast}=1$ then the interior part of the norm $\|\cdot\|_{\ast\ast}$ vanish while for the exterior part we have that
\begin{equation}\label{E35}
\big||z|^{q}D^{q}\tilde{\phi}_{n}(z)\big|\leq1,\ \forall z:=\varepsilon_{n}^{-1}x\in\Omega_{\varepsilon_{n}}\setminus B_{2}.
\end{equation}
Since $\phi_{n}(z)=\tilde{\phi}_{n}(x)$ then $D^{q}\phi_{n}(z)=\varepsilon_{n}^{q}D^{q}\tilde{\phi}_{n}(x)$. Hence by \eqref{E35} we get $$\big||x|^{q}D^{q}\tilde{\phi}_{n}(x)\big|\leq1.$$ If $q=0$ this means that $\|\tilde{\phi}_{n}(x)\|_{L^{\infty}(\Omega\setminus B_{2\varepsilon_{n}})}\leq1$ since $\varepsilon_{n}(\Omega_{\varepsilon_{n}}\setminus B_{2})=\Omega\setminus B_{2\varepsilon_{n}}$. Now, if $q>0$ and $x\in\Omega\setminus B_{\delta'}(0)$ then we have $(\delta')^{q}|D^{q}\tilde{\phi}_{n}(x)|\leq\big||x|^{q}D^{q}\tilde{\phi}_{n}(x)\big|\leq1$. If we replace $q$ with $2q$ we have $\|B_{q}\tilde{\phi}_{n}\|_{L^{\infty}(\Omega\setminus B_{\delta'}(0)}\leq\frac{1}{(\delta')^{2q}}$. Finally, since this is an estimate for each $q$ we have to sum over $q$ and thus we introduce a constant $c_{q}$, dependent only from $q$, in order to have the desired estimate with $c_{\delta'}:=\frac{c_{q}}{(\delta')^{2q}}$. Notice that, since $\mathcal{M}$ is compact, we may always pass to a subsequence, that with slight abuse of notation we continue to call $\xi^{n}:=(\xi_{n}^{i})$, such that $\xi_{n}\rightarrow\xi^{\ast}\in\mathcal{M}$ as $n\rightarrow\infty$.
Then using by estimates \eqref{E34} and by (\cite{axl}, Theorem 2.6) there exists a subsequence still denoted by $\tilde{\phi}_{n}$ such that $\tilde{\phi}_{n}\rightarrow\tilde{\phi}_{0}$ as $n\rightarrow\infty$ in $C^{2m-1,\alpha}$-sense over compact subsets of $\Omega\setminus\{\xi^{\ast}\}$ because $\widetilde{\Omega}\rightarrow\Omega\setminus\{\xi^{\ast}\}$ as $n\rightarrow\infty$. Furthermore, $\tilde{\phi}_{0}$ solves
$$
\left\{ \begin{array}{ll}
         (-\Delta)^{m}\tilde{\phi}_{0}=0& \mbox{in}\ \Omega\setminus\{\xi^{\ast}\},\\
          B_{j}\tilde{\phi}_{0}=0, |j|\leq m-1& \mbox{on}\ \partial\Omega.\
        \end{array} \right.
\vspace{0,2cm}$$
By the first estimate in \eqref{E34} follows that $\|\tilde{\phi}_{0}\|_{L^{\infty}\Omega\setminus\{\xi^{\ast}\}}\leq1$. By Theorem \ref{T2} then the singularity is removable and, thus, $\tilde{\phi}_{0}$ solves
$$
\left\{ \begin{array}{ll}
         (-\Delta)^{m}\tilde{\phi}_{0}=0& \mbox{in}\ \Omega,\\
          B_{j}\tilde{\phi}_{0}=0, |j|\leq m-1& \mbox{on}\ \partial\Omega.\
        \end{array} \right.
\vspace{0,2cm}$$
If $m=1$, $m>1$ with Navier b.c. or $m>1$ even and Dirichlet b.c. integrating by parts we have \underline{directly} that $\tilde{\phi}_{0}=0$. If $1<m$ is odd and with Dirichlet b.c. in order to conclude we need also to use, after an integration by parts, the well-known decomposition of the Laplacian on $\partial\Omega$ valid for any smooth function $u$, namely $$\Delta u=\partial^{2}_{\nu}u+H_{\partial\Omega}\partial_{\nu}u+\Delta_{\tau}u\ \mbox{on}\ \partial\Omega,$$ where $H_{\partial\Omega}\neq0$ denotes the mean curvature at the boundary and $\Delta_{\tau}$ denotes the tangential Laplacian, see (\cite{sperb}, pag.62). For instance, if $m=3$ using by the Dirichlet b.c. the above formula reduces to $\Delta u=\partial^{2}_{\nu}u$ on $\partial\Omega$. Finally, up to a subsequence, $\tilde{\phi}_{n}\rightarrow0$ as $n\rightarrow\infty$ in $C^{2m-1,\alpha}$-sense over compact subsets of $\Omega\setminus\{\xi^{\ast}\}$.
In particular, $$\sum_{|\alpha|\leq2m-1}(\varepsilon_{n}^{-1})^{|\alpha|}|D^{\alpha}\phi_{n}(y)\rightarrow0,\ \mbox{uniformly in}\ |y-(\xi_{n})'|\geq\frac{\delta'}{2\varepsilon_{n}},$$ for any $\delta'>0$ and $i=1,...,k$. We obtain thus that
\begin{equation}\label{E36}
\sum_{i=1}^{k}\sum_{|\alpha|\leq2m-1}\big\|r_{i}^{|\alpha|}D^{\alpha}\phi_{n}\big\|_{L^{\infty}(r_{i}\geq\frac{\delta'}{\varepsilon_{n}})}\rightarrow0, \end{equation}
for any $\delta'>0$. In conclusion, the \emph{exterior part} of the norm $\|\phi_{n}\|_{\ast\ast}$ goes to zero, thus the claim follows.
\\Let us consider now a smooth radial cut-off function $\eta$ with $\eta(s)=1$ if $s<\frac{1}{2}, \eta(s)=0$ if $s\geq1$, and define
$$\hat{\phi}_{n}(y):=\eta\bigg(\frac{\varepsilon_{n}}{\delta_{0}}|y-(\xi_{n}'|\bigg)\phi_{n}(y),$$ such that $$\mbox{supp}\hat{\phi}_{n}\subseteq B_{\frac{\delta_{0}}{\varepsilon_{n}}}(\xi_{n}').$$ Notice that $$\mathcal{L}_{\varepsilon_{n}}(\hat{\phi}_{n})=\eta h_{n}+F(\eta,\phi_{n}),$$ where $$F(f,g):= (-\Delta)^{m}(fg)+C\sum_{i,j=1}^{2m}\frac{\partial^{m}f}{\partial y_{i}\partial y_{j}}\frac{\partial^{m}g}{\partial y_{i}\partial y_{j}}.$$ Thus we get $(-\Delta)^{m}\hat{\phi}_{n}=T_{n}(y)\hat{\phi}_{n}+\eta h_{n}+F(\eta,\phi_{n})$ in $B_{\frac{\delta_{0}}{\varepsilon_{n}}}(\xi_{n}')$ and, since $\hat{\phi}_{n}\in C_{0}^{\infty}(B_{\frac{\delta_{0}}{\varepsilon_{n}}}(\xi_{n}'))$, in particular we have that $B_{j}\hat{\phi}_{n}=0$ on $\partial\Omega$. Since $\|\phi_{n}\|_{\ast\ast}=1$, using \eqref{E36} and the outer estimate in (\cite{clapp}, Lemma 3.2) we have that there exists, up to a subsequence if necessary, an index $i=1,...,k$ such that
\begin{equation}\label{E37}
\liminf_{n\rightarrow\infty}\|\phi_{n}\|_{L^{\infty}(r_{i}<R_{0})}\geq\alpha>0.
\end{equation}
Let us set $\psi_{n,i}(z):=\phi_{n}(z+\xi_{n}')$ where the index $i$ is such that $\sup_{|z-\xi_{n}'|<R_{0}}|\phi_{n}|\geq\alpha>0$. Without loss of generality, we may always assume that the index $i$ is the same for all $n$. Notice that $\psi_{n,i}$ satisfies $$(-\Delta)^{m}\psi_{n,i}-T_{n}(z+\xi_{n}')\psi_{n,i}=h_{n}(z+\xi_{n}')\ \mbox{in}\ \Omega_{n}:=\Omega_{\varepsilon}-\{\xi_{n}'\}.$$
Since $B_{j}\phi_{n}$ are bounded uniformly then as $n\rightarrow\infty$ over compact subsets of $\mathbb{R}^{2m}$ we have $$T_{n}(z+\xi_{n}')=\frac{\alpha_{2m}(2m-1)!}{(1+|z|^{2})^{2m}}(1+o(1)),\ |h_{n}(z+\xi_{n}')|\leq c\|h_{n}\|_{\ast}.$$ Hence, standard elliptic estimates allow us to assume that, as $n\rightarrow\infty, \psi_{n,i}$ converges uniformly over compact subsets of $\mathbb{R}^{2m}$ to a bounded, non-zero solution $\psi$ of $$(-\Delta)^{m}\psi=\frac{\alpha_{2m}(2m-1)!\mu_{i}^{2m}}{(\mu_{i}^{2}+|z|^{2})^{2m}}\psi\ \mbox{in}\ \mathbb{R}^{2m}.$$ This implies that, by Lemma \ref{e}, $\psi$ is a linear combination of the functions $Y_{ij}$ defined in \eqref{E28} and \eqref{E29}.
Orthogonality conditions over $\psi_{n,i}$ pass to the limit thanks to $\|\psi_{n,i}\|_{\infty}\leq1$ and dominated convergence theorem. Thus, this implies that $\psi(y)\equiv0$, a contradiction with \eqref{E37}.
\begin{flushright}
$\Box$
\end{flushright}
Now we will deal with \eqref{E32} lifting the orthogonality constraints $\int_{\Omega_{\varepsilon}}\chi_{i}Z_{i0}\phi=0,\ i=1,...,k$, namely
\begin{equation}\label{E38}
\left\{ \begin{array}{ll}
         \mathcal{L}_{\varepsilon}(\phi)=h& \mbox{in}\ \Omega_{\varepsilon},\\
          B_{j}\phi=0, |j|\leq m-1& \mbox{on}\ \partial\Omega_{\varepsilon},\\
          \int_{\Omega_{\varepsilon}}\chi_{i}Z_{ij}\phi=0& \mbox{for all}\ j=1,...,2m,\ i=1,...,k.
        \end{array} \right.
\vspace{0,2cm}\end{equation}
We have the following a priori estimates for this problem.

\begin{lm}\label{bb}
\end{lm}
\emph{There exist positive constants $\varepsilon_{0}$ and $C$ such that for any solution $\phi$ of problem \eqref{E38} with $h\in L^{\infty}(\Omega_{\varepsilon})$, $\|h\|_{\ast}<\infty$ and with $\xi\in\mathcal{M}$, then for all $\varepsilon\in(0,\varepsilon_{0})$
\begin{equation}\label{E40}
\|\phi\|_{\ast\ast}\leq C|\log\varepsilon|\|h\|_{\ast}.
\end{equation}
}\\$\underline{Proof}$\\
Let $R>R_{0}+1$ be a large, fixed number and $R_{0}$ is the radius of the ball where the cut-off function $\chi$, defined in the previous Lemma, is supported. Let us consider
\begin{equation}
\hat{Z}_{i0}(y):=Z_{0i}(y)-1+a_{i0}G(\varepsilon y,\xi_{i}),
\end{equation}
where $a_{i0}:=\big(H(\xi_{i},\xi_{i})-4m\log(\varepsilon R)\big)^{-1}$. Notice that, by \eqref{E39}, $a_{i0}=\frac{1}{|\log\varepsilon|}$. If $\varepsilon$ is small enough then, by \eqref{E39} and Lagrange Theorem we have $$\hat{Z}_{i0}(y)=Z_{i0}(y)+a_{i0}\bigg(G(\varepsilon y,\xi_{i})-H(\xi_{i},\xi_{i})+4m\log(\varepsilon R)\bigg)$$
\begin{equation}
=Z_{i0}(y)+\frac{1}{|\log\varepsilon|}\bigg(O(\varepsilon r_{i})+4m\log\frac{R}{r_{i}}\bigg),
\end{equation}
and $Z_{i0}(y)=O(1)$. Now, we consider radial smooth cut-off functions $\eta_{1}$ and $\eta_{2}$ with the following properties:
$$0\leq\eta_{1}\leq1,\ \eta_{1}\equiv1\ \mbox{in}\ B_{R}(0),\ \eta_{1}\equiv0\ \mbox{in}\ (B_{R+1}(0))^{c},$$ and $$0\leq\eta_{2}\leq1,\ \eta_{2}\equiv1\ \mbox{in}\ B_{\frac{\delta_{0}}{3\varepsilon}}(0),\ \eta_{2}\equiv0\ \mbox{in}\ \big(B_{\frac{\delta_{0}}{2\varepsilon}}(0)\big)^{c}.$$ Without lack of generality, we may assume $B_{\frac{\delta_{0}}{2\varepsilon}}(0)\subseteq\Omega$. Set $$\eta_{i1}(y):=\eta_{1}(r_{i}),\ \eta_{i2}(y):=\eta_{2}(r_{i}),$$ and define the test function $$\tilde{Z}_{i0}:=\eta_{i1}Z_{i0}+(1-\eta_{i1})\eta_{i2}\hat{Z}_{i0}.$$ Notice that $\mbox{supp}\tilde{Z}_{i0}\subseteq\Omega$ and are all disjoint. Furthermore, intuitively, $\tilde{Z}_{i0}$ resembles the eigenfunction of the operator $\mathcal{L}_{\varepsilon}$ in $\mathbb{R}^{2m}$ with respect to the dilation property of $\mathcal{L}_{\varepsilon}.$ Observe the $\tilde{Z}_{i0}$'s behavior through $\Omega_{\varepsilon}$:
$$
\tilde{Z}_{i0}=
\left\{ \begin{array}{ll}
         Z_{i0}& \mbox{in}\ \Omega_{0}:=\{r_{i}\leq r\},\\
          \eta_{i1}(Z_{i0}-\hat{Z}_{i0})+\hat{Z}_{i0}& \mbox{in}\ \Omega_{1}:=\{R<r_{i}\leq R+1\},\\
          \hat{Z}_{i0}& \mbox{in}\ \Omega_{2}:=\{R+1<r_{i}\leq\frac{\delta_{0}}{3\varepsilon}\},\\
       \eta_{i2}\hat{Z}_{i0}& \mbox{in}\ \Omega_{3}:=\{\frac{\delta_{0}}{3\varepsilon}<r_{i}\leq\frac{\delta_{0}}{2\varepsilon}\},\\
       0& \mbox{otherwise}.
        \end{array} \right.
\vspace{0,2cm}$$
Let $\phi$ be a solution to \eqref{E38}. The main idea of the proof is to modify $\phi$ so that the extra orthogonality conditions with respect to $Z_{i0}$'s hold and try to use Lemma \ref{aa}. Set
\begin{equation}\label{E41}
\hat{\phi}:=\phi+\sum_{i=1}^{k}d_{i}\tilde{Z}_{i0}.
\end{equation}
Our goal is to adjust the constants $d_{i}$ so that
\begin{equation}\label{E42}
\int_{\Omega_{\varepsilon}}\chi_{i}Z_{ij}\hat{\phi}=0,\ \forall j=0,...,2m; i=1,...,k.
\end{equation}
Then,
\begin{equation}\label{E45}
\mathcal{L}_{\varepsilon}(\hat{\phi})=h+\sum_{i=1}^{k}d_{i}\mathcal{L}_{\varepsilon}(\tilde{Z}_{i0})\ \mbox{in}\ \Omega_{\varepsilon}.
\end{equation}
If (\ref{E42}) holds, Lemma \ref{aa} allows us to conclude
\begin{equation}\label{E43}
\|\hat{\phi}\|_{\ast\ast}\leq C\bigg\{ \|h\|_{\ast}+\sum_{i=1}^{k}|d_{i}|\|\mathcal{L}_{\varepsilon}(\tilde{Z}_{i0})\|_{\ast}\bigg \}.
\end{equation}
Estimate (\ref{E40}) is a direct consequence of the following claim.\\
\\$\textbf{Claim}$: \emph{The constant $d_{i}$ are well defined,
\begin{equation}\label{E44}
|d_{i}|\leq C|\log\varepsilon|\|h\|_{\ast}\ \mbox{and}\ \|\mathcal{L}_{\varepsilon}(\tilde{Z}_{0i})\|_{\ast}\leq\frac{C}{|\log\varepsilon|},\ \forall i=1,...,k.
\end{equation}
}
In fact, using by this claim and the fact that $\|\tilde{Z}_{i0}\|_{\ast\ast}\leq C$ we obtain \eqref{E40} as desired. Now, let us prove the claim.
First we find $d_{i}$. From \eqref{E41}, orthogonality conditions \eqref{E42} and the fact that supp$\chi_{j}\eta_{1k}=0$ and supp$\chi_{j}\eta_{2k}=0$ if $i\neq j$ are also satisfied for $\hat{\phi}$ thanks to the fact that $R>R_{0}+1$ we can choose $$d_{i}=-\frac{\int_{\Omega_{\varepsilon}}\chi_{i} Z_{i0}\phi}{\int_{\Omega_{\varepsilon}}\chi_{i} |Z_{i0}|^{2}},\ \forall i=1,...,k$$ and, then, $d_{i}$ is well defined. In order to prove the second inequality in \eqref{E44} we have to compute $\mathcal{L}_{\varepsilon}(\tilde{Z}_{0i})$ in $\Omega_{l},\ l=0,...,3.$ Proceeding exactly as in (\cite{clapp}, Lemma 3.3, Claim 1) we may obtain the desired estimate. Finally, known this, we may prove the first inequality in \eqref{E44}. \\Testing equation \eqref{E45} against $\tilde{Z}_{i0}$ and the above estimate, we get $$|d_{i}|\bigg|\int_{\Omega_{\varepsilon}}\mathcal{L}_{\varepsilon}(\tilde{Z}_{i0})\tilde{Z}_{i0}\bigg|=\bigg|\int_{\Omega_{\varepsilon}}h\tilde{Z}_{i0}+\int_{\Omega_{\varepsilon}}\mathcal{L}_{\varepsilon}(\tilde{Z}_{i0})\hat{\phi}\bigg|$$
$$\leq C\|h\|_{\ast}+C\|\hat{\phi}\|_{\infty}\|\mathcal{L}_{\varepsilon}(\tilde{Z}_{i0})\|_{\ast}.$$ Since $\|\tilde{Z}_{i0})\|_{\ast}=O(1)\ \forall i$ and $\|\hat{\phi}\|_{\infty}\leq\|\hat{\phi}\|_{\ast\ast}$ using by relations (\ref{E43}) and the second inequality in (\ref{E44}) we have \begin{equation}\label{E46}
|d_{i}|\bigg|\int_{\Omega_{\varepsilon}}\mathcal{L}_{\varepsilon}(\tilde{Z}_{i0})\tilde{Z}_{i0}\bigg|\leq C\|h\|_{\ast}+C\sum_{l=1}^{k}\frac{|d_{l}|}{|\log\varepsilon|^{2}}.
\end{equation}
It only remains to estimate the integral term in the left side. Proceeding exactly as in (\cite{clapp}, Lemma 3.3, Claim 2) we may obtain the following claim.\\
\\$\textbf{Claim}$: \emph{If $R$ is large enough, then
\begin{equation}\label{E47}
\bigg|\int_{\Omega_{\varepsilon}}\mathcal{L}_{\varepsilon}(\tilde{Z}_{i0})\tilde{Z}_{i0}\bigg|=\frac{\tilde{C}}{|\log\varepsilon|}(1+o(1)),
\end{equation}
where $\tilde{C}$ is a positive constant independent of $\varepsilon$ and $R$.
}\\
\\At this point, we may replace \eqref{E47} in \eqref{E46} in order to get the desired bounds of $d_{i}$.
\begin{flushright}
$\Box$
\end{flushright}
$\underline{Proof\ of\ Proposition\ \ref{P1}}$\\
First we establish the validity of the a priori estimate \eqref{Eb} for solution $\phi$ of \eqref{E31} with $h\in L^{\infty}(\Omega_{\varepsilon})$ and $\|h\|_{\ast}<\infty$. Lemma \ref{bb} implies
\begin{equation}\label{E48}
\|\phi\|_{\ast\ast}\leq C|\log\varepsilon|\{\|h\|_{\ast}+\sum_{i=1}^{k}\sum_{j=1}^{2m}|c_{ji}|\|\chi_{i}Z_{ji}\|_{\ast}\}.
\end{equation}
Since $\|\chi_{i}Z_{ij}\|_{\ast}\leq C$ then it sufficient to estimate the values of the constant $c_{ij}$. In order to do this, we multiply the first equation in \eqref{E31} by $Z_{ij}\eta_{i2}$, with $\eta_{i2}$ as in Lemma \ref{bb}, and integrating by parts to find
\begin{equation}\label{E50}
\int_{\Omega_{\varepsilon}}\mathcal{L}_{\varepsilon}Z_{ij}\eta_{i2}\phi=\int_{\Omega_{\varepsilon}}hZ_{ij}\eta_{i2}+c_{ij}\int_{\Omega_{\varepsilon}}|Z_{ij}|^{2}\eta_{i2}.
\end{equation}
It is easy to see that $\int_{\Omega_{\varepsilon}}hZ_{ij}\eta_{i2}=O(\|h\|_{\ast}), \int_{\Omega_{\varepsilon}}|Z_{ij}|^{2}\eta_{i2}=C>0$ and $$\bigg|\int_{\Omega_{\varepsilon}}\mathcal{L}_{\varepsilon}Z_{ij}\eta_{i2}\phi\bigg|\leq C\varepsilon|\log\varepsilon|\|\phi\|_{\infty}\leq C\varepsilon|\log\varepsilon|\|\phi\|_{\ast\ast},$$ see also \cite{clapp, delpi,esp}. Using the above estimates in \eqref{E50} we have $$|c_{ij}|\leq C\{\varepsilon|\log\varepsilon|\|\phi\|_{\ast\ast}+\|h\|_{\ast}\},$$ and then we get $$|c_{ij}|\leq C\bigg\{\big(1+\varepsilon|\log\varepsilon|^{2}\big)\|h\|_{\ast}+\varepsilon|\log\varepsilon|^{2}\sum_{l,n}|c_{ln}|\bigg\}.$$ Thus, $|c_{ij}|\leq C\|h\|_{\ast}$ and putting this estimate in \eqref{E48} we conclude.\\ Now, we prove the \emph{solvability assertion}. Notice that Problem \eqref{E38} expressed in a weak form is equivalent to that of finding a $\phi\in\mathcal{H}$, such that $$(\phi,\psi)_{H}=\int_{\Omega_{\varepsilon}}(h+T\phi)\psi,\ \forall \psi\in\mathcal{H},$$ where we consider the Hilbert space $$\mathcal{H}:=\{\phi:\ B_{j}\phi=0\ \mbox{on}\ \partial\Omega_{\varepsilon}\ \mbox{and}\ \int_{\Omega_{\varepsilon}}\chi_{i}Z_{ij}\phi=0,\ \forall i=1,...,k; j=1,...,2m\},$$ endowed with the usual inner product $(\phi,\psi)_{H}$. With the aid of Riesz's representation Theorem, this equation can be rewritten in $\mathcal{H}$ in the operator form $\phi=K(T\phi+h)$, where $K$ is a compact operator in $\mathcal{H}$. Then Fredholm's alternative guarantees unique solvability for any $h$ provided that the homogeneous equation $\phi=K(T\phi)$ has only zero solution in $\mathcal{H}$. This last equation is equivalent to \eqref{E38} with $h\equiv0$. Thus existence of a unique solution follows from the a priori estimate \eqref{E40}. This finishes the proof.
\begin{flushright}
$\Box$
\end{flushright}
\begin{oss}
\end{oss}
\begin{enumerate}
  \item This result implies that the unique solution $\phi=Q(h)$ of \eqref{E31} defines a continuous linear map from the Banach space $C_{\ast}$ of all functions $h\in L^{\infty}(\Omega_{\varepsilon})$ with $\|h\|_{\ast}<+\infty$, into $W^{2m-1,\infty}(\Omega_{\varepsilon})$ with norm uniformly bounded in $\varepsilon$.
  \item The operator $Q$ is differentiable with respect to the variables $\xi'$. In fact, computations similar to those in (\cite{delpi}, pag.17) yield the estimate
 \begin{equation}
\|\partial_{\xi'}Q(h)\|_{\ast\ast}\leq C|\log\varepsilon|^{2}\|h\|_{\ast}.
\end{equation}
This estimate is of crucial importance in the arguments to come.
\end{enumerate}
\section{The intermediate nonlinear problem}\label{S1}
Rather than solve Problem \eqref{E23} directly we shall consider the intermediate nonlinear problem
\begin{equation}\label{E51}
\left\{ \begin{array}{ll}
         \mathcal{L}_{\varepsilon}(\phi)=-R+N(\phi)+\sum_{j=1}^{2m}\sum_{i=1}^{k}c_{ij}\chi_{i}Z_{ij}& \mbox{in}\ \Omega_{\varepsilon},\\
          B_{j}\phi=0, |j|\leq m-1& \mbox{on}\ \partial\Omega_{\varepsilon},\\
          \int_{\Omega_{\varepsilon}}\chi_{i}Z_{ij}\phi=0& \mbox{for all}\ j=1,...,2m,\ i=1,...,k.
        \end{array} \right.
\vspace{0,2cm}\end{equation}
Assuming that conditions in Proposition \ref{P1} hold, we are able to prove the following
\begin{lm}\label{ff}
\end{lm}
\emph{Let $\xi\in\mathcal{M}$. Then, there exists positive constants $\varepsilon_{0}$ and $C$ such that for all $\varepsilon\leq\varepsilon_{0}$ the nonlinear Problem \eqref{E51} has a unique solution $\phi$ which satisfies $$\|\phi\|_{\ast\ast}\leq C\varepsilon|\log\varepsilon|.$$ Moreover, if we consider the map $\xi'\in\mathcal{M}\rightarrow\phi\in C^{2m,\alpha}(\overline{\Omega}_{\varepsilon})$, the derivative $D_{\xi'}\phi$ exists and defines a continuous map of $\xi'$. Besides, there exists a positive constant $\widetilde{C}$ such that $$\|D_{\xi'}\phi\|_{\ast\ast}\leq\widetilde{C}\varepsilon|\log\varepsilon|^{2}.$$}
$\underline{Proof}$\\
In terms of the operator $Q$ defined in Proposition \ref{P1}, \eqref{E51} has the following fixed point representation $$\phi=\mathcal{B}(\phi)\equiv Q(N(\phi)-R).$$ Let us consider the region
$$\mathcal{F}:=\{\phi\in C^{2m,\alpha}(\overline{\Omega}_{\varepsilon}):\ \|\phi\|_{\ast\ast}\leq\varepsilon|\log\varepsilon|\}.$$
From Proposition \ref{P1} we have $$\|\mathcal{B}(\phi)\|_{\ast\ast}\leq C|\log\varepsilon|\{|\ N(\phi)\|_{\ast}+\|R\|_{\ast}\},$$
for arbitrary $\phi$. Hence, by Lemma \ref{b} and Lemma \ref{d} we have $\forall \phi,\phi_{1},\phi_{2}\in\mathcal{F}$
\begin{enumerate}
\item $\|\mathcal{B}(\phi)\|_{\ast\ast}\leq C\varepsilon|\log\varepsilon|;$
\item $\|\mathcal{B}(\phi_{1})-\mathcal{B}(\phi_{2})\|_{\ast\ast}\leq C\varepsilon|\log\varepsilon|^{2}\|\phi_{1}-\phi_{2}\|_{\ast\ast}.$
\end{enumerate}
Then it follows that for all $\varepsilon$ sufficiently small $\mathcal{B}$ is a \emph{contraction mapping} of $\mathcal{F}$, and therefore, by using the Implicit Function Theorem, a unique fixed point of $\mathcal{B}$ exists in this region. The IFT guarantees $C^{1}$ regularity of the map in $\xi'$. Follows exactly the same proof included in (\cite{delpi}, Lemma 4.2) we may show the derivative estimate.
\begin{flushright}
$\Box$
\end{flushright}

\section{Variational reduction}\label{S10}
After Problem \eqref{E51} has been solved, we will find solution to the full Problem \eqref{E23}, or equivalently \eqref{E1}, if we manage to adjust the $k$-uple $\xi'$ in such a way that
\begin{equation}\label{E58}
c_{ij}(\xi')=0,\ \mbox{for\ all}\ i,j.
\end{equation}
A nice feature of this system of equations is that it turns out to be equivalent to finding critical points of a functional of $\xi$ which is close, in appropriate sense, to the energy of the first approximation $W$. Notice that problem \eqref{E58} is indeed variational. In fact, to see that let us consider the energy functional $J_{\rho}$ associated to Problem \eqref{E1}, namely \eqref{E56}. We define the function for $\xi\in\mathcal{M}$
\begin{equation}
\mathfrak{F}_{\varepsilon}(\xi)\equiv J_{\rho}[U(\xi)+\hat{\phi}_{\xi}],
\end{equation}
where $U=U(\xi)$ is our approximate solution from \eqref{E17} and $\hat{\phi}_{\xi}:=\phi(\frac{x}{\varepsilon},\frac{\xi}{\varepsilon}),$ $x\in\Omega$ and $\phi$ the unique solution to Problem \eqref{E51} predicted by Lemma \ref{ff}. Under the assumptions of Lemma \ref{ff}, we obtain, in the following result, that critical points of $\mathfrak{F}_{\varepsilon}$ correspond to solutions of \eqref{E58} for small $\varepsilon$ and, furthermore, the closeness of $\mathfrak{F}_{\varepsilon}$ to $J_{\rho}[U(\xi)]$, for which we know the asymptotic estimate \eqref{E57}.
\begin{lm}\label{gg}
\end{lm}
\emph{The functional $\mathfrak{F}_{\varepsilon}:\mathcal{M}\rightarrow\mathbb{R}$ is of class $C^{1}$. Moreover, for all positive constant $\varepsilon$ small enough, if $D_{\xi}\mathcal{F}_{\varepsilon}(\xi)=0$ then $\xi$ satisfies \eqref{E58}. Besides, for $\xi\in\mathcal{M}$ the following expansion holds
\begin{equation}
\mathfrak{F}_{\varepsilon}(\xi)= J_{\rho}[U(\xi)]+\theta_{\varepsilon}(\xi),
\end{equation}
where $|D_{\xi}^{j}\theta_{\varepsilon}|=o(1)$, uniformly on $\xi\in\mathcal{M}$ as $\varepsilon\rightarrow0, |j|\leq m-1$.}
\\$\underline{Proof}$\\
We define the functional
$$I_{\varepsilon}[w]:=\frac{1}{2}\int_{\Omega_{\varepsilon}}|(-\Delta)^{\frac{m}{2}}w|^{2}dy-\int_{\Omega_{\varepsilon}}V(\varepsilon y)e^{w}dy.$$ Let us differentiate $\mathfrak{F}_{\varepsilon}$ with respect to $\xi$. Notice that, since $J_{\rho}[U(\xi)+\hat{\phi}_{\xi}]=I_{\varepsilon}[W(\xi')+\phi_{\xi'}]$, we can differentiate directly under the integral sign, so that
\begin{eqnarray*}
\partial_{(\xi_{l})_{n}}\mathfrak{F}_{\varepsilon}(\xi)&=&\frac{1}{\varepsilon}DI_{\varepsilon}[W+\phi](\partial_{(\xi_{l})_{n}}W+\partial_{(\xi_{l})_{n}}\phi)\\ &=&\frac{1}{\varepsilon}\sum_{j}\sum_{i}\int_{\Omega_{\varepsilon}}c_{ji}\chi_{i}Z_{ji}(\partial_{(\xi_{l})_{n}}W+\partial_{(\xi_{l})_{n}}\phi).
\end{eqnarray*}
From the results of Section \ref{S1} this expression defines a continuous function of $\xi'$, and hence of $\xi$. Let us assume that $D_{\xi}\mathfrak{F}_{\varepsilon}=0.$ Then, $$\sum_{j}\sum_{i}\int_{\Omega_{\varepsilon}}c_{ij}\chi_{i}Z_{ij}(\partial_{(\xi_{l})_{n}}W+\partial_{(\xi_{l})_{n}}\phi)=0,\ l=1,...,2m, n=1,...,k.$$
Since, by the derivative estimate in Lemma \ref{E51}, we have directly $$(\partial_{(\xi_{l})_{n}}W+\partial_{(\xi_{l})_{n}}\phi)Z_{ln}+o(1),$$ where $o(1)$ is uniformly small as $\varepsilon\rightarrow0$. Thus, we get that $D_{\xi}\mathfrak{F}_{\varepsilon}=0$ implies the validity of the equations in term of $\ast\ast$-norm
$$\sum_{j}\sum_{i}\int_{\Omega_{\varepsilon}}c_{ij}\chi_{i}Z_{ij}(Z_{ln}+o(1))=0.$$ This system is dominant diagonal, thus we get $c_{ij}=0$ for all $i,j$. \\The closeness property follows directly from an application of Taylor expansion for $\mathfrak{F}_{\varepsilon}$ in $\Omega_{\varepsilon}$ and from the estimates in Lemma \ref{E51} See also (\cite{delpi}, Lemma 5.2).
This concludes the proof.

\begin{flushright}
$\Box$
\end{flushright}
\section{Proof of Theorem \ref{T1}}\label{S20}
Taking into account Lemma \ref{gg}, a solution to \eqref{E1} exists if we prove the existence of a critical point of $\mathfrak{F}_{\varepsilon}$, which automatically implies that $c_{ij}=0$ for all $i,j$. Next, the qualitative properties of the solution found follow from the chosen ansatz.
Finding critical points of $\mathfrak{F}_{\varepsilon}(\xi)$ is equivalent to finding critical points of
\begin{equation}
\widetilde{\mathfrak{F}}_{\varepsilon}:=\mathfrak{F}_{\varepsilon}-4mb_{m}k|\log\varepsilon|.
\end{equation}
On the other hand, if $\xi\in\mathcal{M}$, from Lemma \ref{c} and Lemma \ref{gg} we get the existence of universal constants $\alpha>0$ and $\beta$ such that
\begin{equation}\label{E60}
\alpha\widetilde{\mathfrak{F}}_{\varepsilon}+\beta=\varphi_{k}(\xi)+O(\varepsilon).
\end{equation}
With the same argument in (\cite{clapp}, Theorem 2) and (\cite{delpi}, Theorem 2), we may prove that, under the assumptions of Theorem \ref{T1}, $\widetilde{\mathfrak{F}}_{\varepsilon}$ has a critical point in $\mathcal{M}$ for $\varepsilon$ small enough. By \eqref{E60} the proof is concluded.
\begin{flushright}
$\Box$
\end{flushright}
\begin{center}
\textbf{\large{Appendix}}
\end{center}
In the sequel we calculate the explicit values $c_{0}$ and $c_{1}$ used in Lemma \ref{c}. This values, appeared until now only in the case $m=1,2$, are broadly used when we work with elliptic equations with exponential nonlinearity.
\begin{lm}
\end{lm}
\emph{\begin{enumerate}
  \item $\int_{\mathbb{R}^{2m}}\frac{dy}{(1+|y|^{2})^{2m}}=\frac{\pi^{m}(m-1)!}{(2m-1)!}$,\\
  \item $\int_{\mathbb{R}^{2m}}\frac{\log(1+|y|^{2})dy}{(1+|y|^{2})^{2m}}=\frac{\pi^{m}(m-1)!}{m(2m-1)!}$.
\end{enumerate}}
$\underline{Proof}$\\
We start to calculate the following integral $\int_{\mathbb{R}^{2m}}\frac{dy}{(1+|y|^{2})^{\alpha}}\ \mbox{for}\ \alpha>m$, interesting by itself.
The formula of integration in spherical coordinates of an integrable rotationally symmetric function yields
\begin{displaymath}
\int_{\mathbb{R}^{2m}}\frac{dy}{(1+|y|^{2})^{\alpha}}=\int_{0}^{\infty} \bigg(\int_{\partial B(0,r)}\frac{d\sigma}{(1+r^{2})^{\alpha}}\bigg)dr.
\end{displaymath}
We observe that
\begin{displaymath}
\int_{\partial B(0,r)}\frac{d\sigma}{(1+r^{2})^{\alpha}} = \frac{1}{(1+r^{2})^{\alpha}} 2m\omega_{2m}r^{2m-1}= \frac{2\pi^{m}r^{2m-1}}{\Gamma(m)(1+r^{2})^{\alpha}}.
\end{displaymath}
By the change of variables $r=\tan\theta$,
\begin{displaymath}
\int_{0}^{\infty}\frac{r^{2m-1}}{(1+r^{2})^{\alpha}}dr=\int_{0}^{\frac{\pi}{2}}(\cos\theta)^{2(\alpha-m)-1}(\sin\theta)^{2m-1}d\theta.
\end{displaymath}
Recalling the definition of the \emph{Beta function} $\beta(s,t):=2\int_{0}^{\frac{\pi}{2}}(\cos\theta)^{2s-1}(\sin\theta)^{2t-1}d\theta$ for all $(s,t)\in{Q}^{+}$, the set of positive rational numbers, we have
\begin{displaymath}
\int_{0}^{\infty}\frac{r^{2m-1}}{(1+r^{2})^{\alpha}}dr=\frac{1}{2}\beta(\alpha-m,m).
\end{displaymath}
Recalling the relation between Beta and Gamma function, $\beta(s,t)=\frac{\Gamma(s)\Gamma(t)}{\Gamma(s+t)}$ we have
$
\int_{0}^{\infty}\frac{r^{2m-1}}{(1+r^{2})^{\alpha}}dr=\frac{1}{2}\frac{\Gamma(\alpha-m)\Gamma(m)}{\Gamma(\alpha)}$
and thus,
\begin{displaymath}
\int_{\mathbb{R}^{2m}}\frac{dy}{(1+|y|^{2})^{\alpha}}= \frac{\pi^{m}\Gamma(\alpha-m)}{\Gamma(\alpha)}.\\
\end{displaymath}
In particular, with $\alpha=2m$  and by \emph{Legendre duplication formula} for $\Gamma$:
\begin{displaymath}
\Gamma(m)\Gamma\bigg(m+\frac{1}{2}\bigg) = 2^{1-2m}\sqrt{\pi}\Gamma(m)
\end{displaymath}
we have
\begin{displaymath}
\int_{\mathbb{R}^{2m}}\frac{dy}{(1+|y|^{2})^{2m}}= \frac{\pi^{m}\Gamma(m)}{\Gamma(2m)} = \frac{2\pi^{m+\frac{1}{2}}}{2^{2m}\Gamma(m+\frac{1}{2})}.
\end{displaymath}
Using the classical equality $\Gamma(m+\frac{1}{2}) = \frac{\sqrt{\pi}(2m-1)!!}{2^{m}}$, the fact that $(2m-1)!! = \frac{(2m)!}{m!2^{m}}$ and by the definition of the factorial of a nonnegative integer the proof of the first point is concluded.\\ Moreover, integrating by parts, we note that
\begin{displaymath}
\int_{\mathbb{R}^{2m}}\frac{\log(1+|y|^{2})dy}{(1+|y|^{2})^{2m}}=\frac{2m\omega_{2m}}{m}\int_{\mathbb{R}^{2m}}\frac{dy}{(1+|y|^{2})^{2m}}=2\omega_{2m}\int_{\mathbb{R}^{2m}}\frac{dy}{(1+|y|^{2})^{2m}}
\end{displaymath}
and the proof is concluded.
\begin{flushright}
$\Box$
\end{flushright}

\bibliographystyle{plain}

\end{document}